\documentclass[reqno,12pt]{amsart}

\usepackage[utf8]{inputenc}

\usepackage{enumerate}
\usepackage[margin=1in,marginparwidth=0.7in]{geometry}
\usepackage{ifpdf}
\usepackage{amsmath}
\usepackage{amsfonts}
\usepackage{amssymb}
\usepackage{amsthm}
\usepackage[ocgcolorlinks,hyperfootnotes=false,colorlinks=true,citecolor=blue,linkcolor=blue,urlcolor=blue]{hyperref}
\usepackage{setspace}
\usepackage{amsrefs}
\usepackage{nicefrac}
\usepackage{graphicx}
\usepackage{xcolor}
\usepackage{mathtools}

\newcommand{\ignore}[1]{}

\renewcommand{\Re}{\operatorname{Re}}
\renewcommand{\Im}{\operatorname{Im}}

\newcommand{\supp}{\operatorname{supp}}

\newcommand{\sabs}[1]{\lvert {#1} \rvert}

\newcommand{\snorm}[1]{\lVert {#1} \rVert}

\newcommand{\C}{{\mathbb{C}}}
\newcommand{\R}{{\mathbb{R}}}

\newcommand{\N}{{\mathbb{N}}}

\newcommand{\D}{{\mathbb{D}}}

\newcommand{\bT}{{\mathbb{T}}}

\newcommand{\sF}{{\mathcal{F}}}

\newcommand{\sO}{{\mathcal{O}}}

\newcommand{\sY}{{\mathcal{Y}}}

\newcommand{\CR}{{\mathit{CR}}}
\renewcommand{\DH}{{\operatorname{DH}}}

\newcommand{\SADH}{{\operatorname{SADH}}}

\newtheorem{thm}{Theorem}[section]

\newtheorem{cor}[thm]{Corollary}

\newtheorem{lemma}[thm]{Lemma}

\theoremstyle{definition}
\newtheorem{defn}[thm]{Definition}
\newtheorem{example}[thm]{Example}

\theoremstyle{remark}
\newtheorem{remark}[thm]{Remark}

\newcommand{\avoidbreak}{\postdisplaypenalty=100}

\author{Ji\v{r}\'{\i} Lebl}
\thanks{The first author was in part supported by Simons Foundation collaboration grant 710294.}
\address{Department of Mathematics, Oklahoma State University,
Stillwater, OK 74078, USA}
\email{lebl@okstate.edu}

\author{Alan Noell}
\address{Department of Mathematics, Oklahoma State University,
Stillwater, OK 74078, USA}
\email{alan.noell@okstate.edu}

\author{Sivaguru Ravisankar}
\thanks{The last author was in part supported by Mathematical Research Impact Centric Support (MATRICS) grant MTR/2022/000865 from the Science and Engineering Research Board (SERB), Government of India.}
\address{Tata Institute of Fundamental Research, Centre for Applicable Mathematics, Bengaluru 560065, India}
\email{sivaguru@tifrbng.res.in}

\date{\today}

\title{CR functions at CR singularities: approximation, extension, and hulls}


\begin{document}

\begin{abstract}
We study three possible definitions of the notion of CR functions at CR singular
points, their extension to a fixed-neighborhood of
the singular point, and analogues of the
Baouendi--Tr\`eves approximation in a fixed neighborhood.
In particular,
we give a construction of certain disc hulls,
which, if large enough, give
the fixed-neighborhood extension and approximation
properties.
We provide many examples showing the distinctions between the
classes and the various properties studied.
\end{abstract}

\maketitle



\section{Introduction} \label{section:intro}

Let $M \subset \C^n$ be a real submanifold  and
$T_\eta^{0,1}M$ the span of the antiholomorphic vectors at a point $\eta \in
M$.  The manifold $M$ is said to be CR at $q$ if the dimension of 
$T_\eta^{0,1}M$ is constant near $q$.  A natural generalization of
holomorphic functions is a so-called CR function, a function killed
by $T^{0,1}M$ vector fields.  If $M$ and $f$ are CR and real-analytic,
then $f$ extends locally to a neighborhood as a holomorphic function
by the theorem of Severi~\cite{Severi:31}.
If the regularity is lower, we do not always get such an extension, but
basic questions of when and into what set extension holds are relatively well-understood: see the pioneering work by Lewy~\cite{Lewy}, Kohn--Rossi~\cite{Kohn-Rossi}, Tumanov~\cite{Tumanov}, and many others.
The problem then is to understand the CR singular setting.  In particular,
it is not immediately clear what is the most natural notion of
\emph{CR function}.  One possible definition is simply to consider
functions that are CR at CR points, or equivalently, those that are killed
by vector fields valued in $T_\eta^{0,1}M$ at each point.  We call these
CR functions, and the authors have studied their extension properties
in \cites{crext1,crext2,LNR:Severi}.

A related question is the approximation of functions by polynomials or
entire functions.  In the totally real case, see
H\"{o}rmander--Wermer \cite{HormanderWermer},
Harvey--Wells \cites{HarveyWells71,HarveyWells72},
and in general
the survey article by Dwilewicz \cite{Dwilewicz}.
The celebrated Baouendi--Tr\`eves theorem \cite{B-T} applied to CR functions on CR
manifolds says that CR functions can be approximated in a fixed neighborhood
(not depending on the function)
of any point by holomorphic polynomials. At CR singular points, such a theorem does not necessarily hold for the CR functions as we defined them above.  We therefore define the class of CR${}_P$ functions to be those that are locally uniform limits of holomorphic polynomials. 
Mondal~\cite{Mondal}, extending the work of Mergelyan~\cite{Mergelyan} and
Minsker~\cite{Minsker}, recently studied the
approximation property for continuous functions on certain
CR singular manifolds that are totally real at CR points.
Finally, we write CR${}_H$ for functions that are restrictions of
holomorphic functions, that is, those that do extend to some neighborhood.
Sometimes $\sO$ is used for such functions, but we will reserve
$\sO(U)$ for functions holomorphic on a specific open set $U$ for clarity.
See Section~\ref{section:crfuncs} for precise definitions and the statement
of one of our main results.


There are other possible definitions for what constitutes CR functions
on a CR singular manifold.  For example, 
Nacinovich-Porten~\cite{NacPor} define a class of functions by considering
the local closure of germs of CR${}_H$ functions, and they study its extension
properties.
Their definition is equivalent to
iterating (possibly infinitely many times) an approximation procedure,
and hence this class lies in between our CR${}_H$ and CR${}_P$.

A subtle but important issue with respect to extension and approximation
of CR functions is the size of the neighborhood, namely, whether the
neighborhood to which the function extends, or on which it is approximated,
depends on the function itself (or not, as is the case for Baouendi--Tr\`eves).
For this purpose, we define four different properties a class of functions can
satisfy:
the \emph{extension property} and \emph{approximation property}
for functions that extend to be holomorphic to a neighborhood and those
that are approximable locally uniformly by holomorphic polynomials,
and the \emph{fixed-neighborhood extension property} and
the \emph{fixed-neighborhood approximation property} if extension
or approximation can be done in a neighborhood independent of the
function. The extension property implies the approximation property, but not
vice-versa. See Section~\ref{section:extapproxhulls} for the definitions and basic examples.

One can use families of discs to extend CR${}_P$ functions.
If we can fill a neighborhood of a point with analytic
discs attached to the submanifold, possibly iterating this
construction, we can extend CR${}_P$ functions to this neighborhood.
A subtle issue is that we need these discs to be attached to an
arbitrarily small neighborhood of a point, that is, we do not
a priori have the fixed-neighborhood approximation property
for CR${}_P$ functions at a CR singular point.
We formalize this property in Definition~\ref{defn:DHcondition}
and define what we call the \emph{DH condition}.

If we can further shrink such discs to a point, we can
apply the Kontinuit\"atssatz
(see Ivashkovich \cite{Ivashkovich}) to get analytic continuation
of CR${}_H$ functions, that is, functions holomorphic in
some neighborhood, and we can hope to get a fixed-neighborhood
extension property.  However, showing that these continuations
are single-valued is difficult.  For points in the 
\emph{shrinking approximate disc hull}, or $\SADH_q$,
we can shrink these discs nicely to the point $q$.
Here we require only approximately attached discs.
If $\SADH_q$ contains a neighborhood of the point $q$,
then we get the fixed-neighborhood extension property
(see Corollary~\ref{cor:bigSADHnoiter}).  However, if
one requires iteration, much more care has to be taken.
We must require the resulting paths in the iteration not to
intersect, so we define the \emph{iterated SADH condition}
in Definition~\ref{defn:SADHcondition}.  With this condition,
we obtain the fixed-neighborhood extension property in
Theorem~\ref{thm:prophullSADH}.

These ideas can be
combined in many different ways with existing results such as
the theorem of Hanges and Tr\`eves on the propagation of extension
along complex analytic curves through CR points.
See Section~\ref{section:results} for these results.
In the CR singular case,
a large disc hull and in fact
the strong iterated SADH condition
can appear even
in the Levi-flat case,
where near the CR points all discs lie in the manifold.
In Section~\ref{section:Hull}, we show that
the manifold given by $w = \bar{z}_1 z_2$, which is Levi-flat at CR points
(in fact, an image of $\R^2 \times \C$), nevertheless satisfies
the strong iterated SADH condition at the origin.

We are interested in producing examples showing that the given
classes of functions are distinct. In particular, we wish to
extend the Baouendi--Tr\`eves result to a more general CR singular setting.
A natural question is whether a class of CR functions on
a submanifold of $\C^n \times \R$
has the approximation property.  In this setting, one may guess that
the independence of the holomorphic and real coordinates may be used
in tandem and that an approximation theorem holds as long as it holds whenever
the last coordinate is fixed.  That is, perhaps one can combine the classical Weierstrass theorem with Mergelyan's theorem.
It turns out (see Theorem~\ref{thm:graph}) that such a result holds with an additional hypothesis.
This result shows that fixed-neighborhood extension is not necessary for fixed-neighborhood
approximation (see, for example, 
Section~\ref{section:BTnofixednbhd}).
With such a result, one can prove approximation results for flat hyperbolic
Bishop surfaces (see Section~\ref{section:HyperbolicBishop}).
For flat elliptic Bishop surfaces, in Section~\ref{section:BishopFlatElliptic}
we prove approximation under the extra assumption of extension to the
natural family of attached analytic discs.
For the special elliptic Bishop surface $w = \sabs{z}^2$,
we can adapt the proof of Baouendi--Tr\`eves
for the approximation and obtain a linear operator via
integration. 
See Section~\ref{section:BTBishop}.

%





\section{CR functions} \label{section:crfuncs}

Let $M \subset \C^n$ be a real submanifold of regularity at least $C^1$.
We will assume that all submanifolds are embedded.
A point $q \in M$ is called a \emph{CR point}
if the dimension of
\begin{equation}
T_\eta^{0,1}M = \C \otimes T_\eta M \cap \operatorname{span}_{\C} \left
\{
\frac{\partial}{\partial\bar{z}_1}\Big|_{\eta},
\ldots ,
\frac{\partial}{\partial\bar{z}_n}\Big|_{\eta}
\right\}
\end{equation}
is constant as a function of $\eta$ in some neighborhood of $q$.
Write $M_{\CR} \subset M$
for the set of CR points.  A point $q \in M$ is said to be a \emph{CR singular point}
if $q \notin M_{\CR}$.  A submanifold is said to be CR if it has no CR
singular points, or in other words if $M = M_{\CR}$.
A vector field $L \in \Gamma(\C \otimes TM)$ is said to be a \emph{CR vector field}
if $L_q \in T_q^{0,1} M$ for all $q \in M$.
We remark that our definition of CR vector field
includes vector fields at CR singular points
as well.

There are several natural definitions of what it means for a function to be CR on a possibly
CR singular submanifold.  That is, there are different ways to define the analogue of
holomorphic functions on a real submanifold.  We focus on three
such definitions.

\begin{defn}
Suppose $M \subset \C^n$ is a
real $C^1$ submanifold and $f \colon M \to \C$ is a continuous function. We say:
\begin{enumerate}[(i)]
\item
$f$ is a \emph{CR function} if $Lf = 0$ for every (continuous) CR vector field $L$ on $M$, interpreted in terms of distributions if $f$ is only continuous.
\item
$f$ is a \emph{CR${}_P$ function} if for every $q \in M$ there exist a
compact neighborhood $K \subset M$ of $q$ and a sequence of holomorphic polynomials
$\{ P_j \}$ such that $P_j$ converges uniformly on $K$ to $f|_K$.
\item
$f$ is a \emph{CR${}_H$ function} if for every $q \in M$ there exist a
neighborhood $U \subset \C^n$ of $q$ and a holomorphic function $F \colon
U \to \C$ such that $F|_{M \cap U} = f|_{M \cap U}$.
\end{enumerate}

We write $\CR^k(M)$ for the set of CR functions that are in $C^k(M)$, including $k=0$ for continuous, $k=\infty$ for smooth, and $k=\omega$ for real-analytic.
We define $\CR^k_P(M)$ in a similar way.  As CR${}_H$ functions are always
of the same regularity as the manifold, we will write simply
$\CR_H(M)$.
\end{defn}

\begin{remark}
The class $\CR_H(M)$ makes sense for any set $M$, not necessarily a
manifold.  Correspondingly, some of the basic results in this paper hold without $M$
necessarily being a manifold.
Similarly the classes $\CR^k_P(M)$ also make
sense for a non-manifold $M$, although the regularity of the functions then
needs to be understood, for example, in the sense of Whitney.
We leave such generalizations to the reader.
\end{remark}

It is easy to see that, for a $C^1$ submanifold $M$, the
set of CR points $M_{\CR}$ is an open dense set
in $M$. Therefore,
a function $f$ is CR if and only if
$f|_{M_{\CR}}$ is a CR function on $M_{\CR}$.
We remark that if $f$ is a CR${}_P$ function, then
$Lf = 0$ in the sense of distributions for every
CR vector field $L$, and hence $f$ is CR.
So CR${}_P$ implies CR for all $C^1$ submanifolds,
but for CR singular submanifolds, the converse
may not be true (see below). Note that
CR$_H$ trivially implies CR. In fact, CR$_H$ implies CR${}_P$ because we can use a series expansion at each point.

When $M$ is real-analytic and CR, by Severi's theorem \cite{Severi:31} real-analytic CR functions are
restrictions of holomorphic functions.  That is, in this
case
$\CR^\omega(M) = 
\CR^\omega_P(M) = 
\CR_H(M)$.
On the other hand, there do exist smooth CR functions on CR submanifolds
that are not restrictions of holomorphic functions.
In the presence of CR singularities, these classes can be distinct even in
the real-analytic case.

For CR submanifolds of class $C^2$, the first two definitions are equivalent, which follows
from the Baouendi--Tr\`eves approximation
theorem \cite{B-T}: \emph{If $M$ is a CR submanifold of class $C^2$ and $q \in M$,
then there exists a compact
neighborhood $K \subset M$ of $q$ such that for every CR function $f$ on $M$ there exists a sequence 
$\{ P_j \}$ of holomorphic polynomials converging uniformly on $K$ to $f|_K$.}
So for a CR submanifold of class $C^2$, $\CR^k(M)=\CR^k_P(M)$ for all $k$.

A key point in the Baouendi--Tr\`eves theorem is that the neighborhood $K$
is independent of $f$; it depends only on $M$ and the point $q$.  We will
see that, in the CR singular case, there exist $M$ for which the
conclusion of the Baouendi--Tr\`eves theorem does not hold even for
CR${}_P$ functions. (See Theorem \ref{thm:elliptic} or Theorem \ref{thm:anote}.)


When considering regularity of manifolds or functions,
we use the order $0 < 1 < 2 < \cdots < \infty < \omega$.

\begin{thm}\label{thm:classes}
Let $M \subset \C^n$ be a real submanifold of regularity $C^\ell$
for $\ell \geq 1$.
\begin{enumerate}[(i)]
\item
$\CR^k(M) \supset \CR^k_P(M) \supset \CR_H(M)$ for all $k \leq \ell$.
\item
There exists a real-analytic submanifold $M$ such that,
for every $k$, $\CR^k(M) \supsetneq \CR^k_P(M)$. 
\item
There exists a real-analytic submanifold $M$ such that
$\CR^\omega_P(M) \supsetneq \CR_H(M)$.
\end{enumerate}
\end{thm}

\begin{proof}
(i) This follows from the earlier observations that CR${}_P$ implies CR and CR$_H$ implies CR${}_P$.

(ii) See Theorem \ref{thm:elliptic} or Theorem \ref{thm:tar}.

(iii) See Theorem \ref{thm:hyperpara}.
\end{proof}

The space $\CR_H(M)$ is defined via local extension, but sometimes it is necessary to have one global extension. Recall that a real submanifold is \emph{generic} at a point if the complex
differentials of its defining functions are linearly independent over $\C$ at that point.

\begin{lemma} \label{lemma:simplyconn}
Suppose $M \subset \C^n$ is a $C^1$ real submanifold
that is either generic at every CR point or
simply connected.
Then for every $f \in \CR_H(M)$ there exist an open neighborhood
$U$ of $M$ in $\C^n$ and a holomorphic function $F$ on $U$ such that $F|_M = f$.
\end{lemma}

\begin{proof}
If $M$ is generic at every CR point, then the extension is unique
locally at each CR point; as those points are dense, the extension is unique everywhere, and the result follows.
If $M$ is simply connected, the result follows by the monodromy theorem.
\end{proof}


\section{Extension, approximation, and hulls} \label{section:extapproxhulls}

\begin{defn}
Suppose $M \subset \C^n$ is a real submanifold, $q \in M$,
and $\sF$ is a class of functions on $M$.  We say:
\begin{enumerate}[(i)]
\item
$M$ has the \emph{extension property for $\sF$ at $q$}
if for every $f \in \sF$ there exist a neighborhood $U$ of $q$
in $\C^n$ and a holomorphic function $F \colon U \to \C$
such that $F|_{U \cap M} = f|_{U \cap M}$.
\item
$M$ has the \emph{fixed-neighborhood extension
property for $\sF$ at $q$}
if there exists a neighborhood $U$ of $q$ in $\C^n$
such that
for every $f \in \sF$
there exists a holomorphic function $F \colon U \to \C$
such that $F|_{U \cap M} = f|_{U \cap M}$.
\item
$M$ has the \emph{approximation property
for $\sF$ at $q$} if for every $f \in \sF$
there exists a compact neighborhood $K$ of $q$
in $M$ such that $f$ is the uniform limit
on $K$ of a sequence of holomorphic polynomials.
\item
$M$ has the \emph{fixed-neighborhood
approximation property
for $\sF$ at $q$} if the following analogue of the
Baouendi--Tr\`eves approximation theorem holds at $q$ for
functions in $\sF$:
There exists a compact neighborhood $K$ of $q$
in $M$ such that every $f \in \sF$ is the uniform limit
on $K$ of a sequence of holomorphic polynomials.
\end{enumerate}
When we say simply that $M$ has
one of the properties above without mentioning
a point $q$, we mean it has the property at all points. If the submanifold is given, we may say that the class $\sF$ has the indicated property.
\end{defn}

We note that Nacinovich-Porten \cite{NacPor} have studied the extension
and approximation properties for a class of functions in between CR${}_H$ and CR${}_P$.

We make some immediate observations:
$\CR_H(M)$ always has the extension property and $\CR^k_P(M)$ always has the approximation property.
The fixed-neighborhood extension property for a class at a point implies
the fixed-neighborhood approximation property for that class at that point.
If $\CR^k(M)$ has the fixed-neighborhood approximation property,
then $\CR^k(M)=\CR^k_P(M)$.
These properties are invariant under holomorphic changes of coordinates.
If $M$ is contained in the Levi-flat hypersurface given by $\Im z_n = 0$,
then $\CR_H(M)$ does not have the fixed-neighborhood extension property as
$1/(z_n-i\epsilon)$ is in $\CR_H(M)$ for all real $\epsilon \not= 0$.

Note that the fixed-neighborhood approximation property for $\CR_H(M)$ need not imply the
fixed-neighborhood extension property for $\CR_H(M)$.
See any one of Theorems \ref{thm:elliptic}, \ref{thm:hyperpara}, and \ref{thm:anote}, or the CR case.


A standard procedure (although it is not sufficient) to construct the polynomial hull is to consider the so-called disc hull.  For some sets,
we may also have to iterate this procedure, as Example~\ref{example:neediterated} shows.
Some of the following definitions and examples are stated for an arbitrary
set $X$ rather than a submanifold.  The reason is that if we iterate the
given constructions, we will obtain sets that are not necessarily
submanifolds in the intermediate steps.

Let $\D \subset \C$ denote the unit disc.
By an \emph{analytic disc attached to $X \subset \C^n$} we mean
a continuous function
$\varphi \colon \overline{\D} \to \C^n$ that is holomorphic
on $\D$ and satisfies $\varphi(\partial \D) \subset X$.
If we say $\varphi$ is an analytic disc \emph{through $p$} we mean in addition that
$p \in \varphi(\D)$. 

\begin{defn} \label{defn:DHcondition}
Let $X \subset \C^n$.  Define
\begin{align}
\DH(X) = &
\{ z \in \C^n :
\exists\text{ an analytic disc attached to $X$ through $z$} \}, \\
\DH^k(X) = &
\underbrace{\DH(\cdots \DH(\DH(X))\cdots)}_{k \text{ times}}
\end{align}
We call $\DH(X)$ the \emph{disc hull of $X$}
and
$\DH^k(X)$ the \emph{$k$-fold iterated disc hull of $X$}.

We say that $X$ satisfies the \emph{DH condition at $q \in X$} if for every
neighborhood $U \subset X$ of $q$, there is a $k$ such that $\DH^k(U)$
is a (not necessarily open) neighborhood of $q$ in $\C^n$.
\end{defn}

The set $\DH^k(X)$ is a subset of the polynomial hull.
We wish to
apply the Kontinuit\"atssatz to functions defined on a neighborhood of $X$,
and for this purpose we also need
to be able to continuously shrink these discs.
It is not always possible to shrink
the discs that make up the disc hull even if $X$ is a submanifold
(see the examples below),
so we require another definition.
On the other hand, we do not need to exactly attach these discs, as we only
need their boundaries to be near $X$ to apply the Kontinuit\"atssatz
to functions in $\CR_H$. Thus, we can weaken the
attachment and only consider approximately attached discs.

\begin{defn}  \label{defn:SADH}
Let $q \in X \subset \C^n$.  For $\epsilon>0$ let $X_\epsilon$ denote the
$\epsilon$-neighborhood of $X$.  Define
\begin{align}
\SADH_q(X) = 
\{ z \in \C^n : {} & \text{for each } \epsilon > 0,
\exists\text{ a continuous family of analytic discs } \varphi_t \colon \overline{\D} \to \C^n,  \notag
\\
&
t \in [0,1] , ~
z = \varphi_1(0) , ~
\varphi_t(\partial \D) \subset X_{\epsilon} \; \forall t \in [0,1] , ~
\varphi_0 \equiv q , ~ \text{and} \notag
\\
&
\snorm{\varphi_t(0)-q} \text{ is a strictly increasing function of $t$} \},
\end{align}
We call $\SADH_q(X)$ the
\emph{shrinking approximate disc hull of $X$ at $q$}.
\end{defn}


We remark that the content of Definition \ref{defn:SADH} is unchanged
if we require $z \in \varphi_1(\D)$ instead of $\varphi_1(0)=z$.

Finally, we need to take into account single-valuedness of the extension,
especially if we are planning to iterate this construction.  It is possible to
extend functions in $\CR_H(X)$ to a ball in $\SADH_q(X)$ if it in fact
contains such a ball, but it is not possible simply to iterate
this procedure. That is, just because 
$\SADH_q(\SADH_q(X))$ contains a ball does not mean we can extend
to this ball---we need extra assumptions.
One possible
way to do this is to ensure that the paths along which we can extend
do not create loops.  With an eye towards the application
to weighted homogeneous submanifolds, and specifically the
example in Theorem~\ref{thm:tar},
we make the following definition, which will allow iterating
the shrinking disc hulls.

\begin{defn} \label{defn:SADHcondition}
We say $Y$ is a \emph{union of a compact family of nonintersecting paths from
$\SADH_q(X)$} if there is a family $\sY$ of paths $\psi(t) = \varphi_t(0)$
arising from the definition of $\SADH_q(X)$, $\sY$ is compact with respect
to the uniform norm, $Y$ is the union of the images of paths in $\sY$, and
whenever
$\psi_1,\psi_2 \in \sY$ are such that $\psi_1(t_1)=\psi_2(t_2)$ for some
$t_1,t_2 > 0$, then $\psi_1([0,t_1]) = \psi_2([0,t_2])$.

We will call the sets
$X_0, X_1, \ldots, X_k$ 
\emph{SADH iterates at $q$}
if $X_0 \subset X$ is a neighborhood of $q$,
and each $X_j$ for $j=1,\ldots,k$
is a union of a compact family of nonintersecting paths from
$\SADH_q(X_{j-1})$.

We say that $X$ satisfies the \emph{iterated SADH condition at $q$}
if $X_0 = X$
and there are sets $X_1, X_2, \ldots, X_k$ that are SADH iterates at $q$
and $X_k$ contains a neighborhood of $q$ in $\C^n$.

We say that $X$ satisfies the \emph{strong iterated SADH condition at $q$}
if $X$ satisfies the iterated SADH condition at $q$ for all
neighborhoods $X_0 \subset X$ of $q$.
%
\end{defn}

We will show that if $\SADH_q(X)$ contains a ball centered at $q$,
then functions from $\CR_H(X)$ will extend to that ball.
However,
it is not true that functions in $\CR_H(X)$ necessarily extend
holomorphically to a neighborhood of $\SADH_q(X)$.  We can use
the Kontinuit\"atssatz, but that only gives us analytic continuation
along the paths arising from the disc families, not unrestricted continuation.
Without some condition on $\SADH_q(X)$ (or a subset), it is not always
possible to have functions in $\CR_H(X)$ to extend to
$\CR_H\bigl(\SADH_q(X)\bigr)$ precisely because we do not get unrestricted
continuation.  See Example~\ref{example:nonuiniqueonSADH}.

A few remarks are in order.

\begin{remark}
The metric used in the definition of $\SADH_q(X)$ need not be the
Euclidean distance; other metrics could easily be used.
\end{remark}

\begin{remark}
The way we defined $\SADH_q(X)$ used paths of strictly increasing distance.
We get the same set if we relax this definition simply to increasing distance:
If $p$ can be reached by a family of discs where $\snorm{\varphi_t(0)-q}$
is only increasing and with boundary in an $(\epsilon/2)$-neighborhood of $X$,
we can perturb the family to get a strictly increasing 
$\snorm{\varphi_t(0)-q}$ and boundary in an $\epsilon$-neighborhood of $X$.
Assuming $\snorm{\varphi_t(0)-q}$ is strictly increasing makes some of
the arguments somewhat cleaner, and as we just saw, there is no loss in generality.
\end{remark}

We will generally apply $\SADH_q$ to neighborhoods of $q$ in a manifold. A manifold is always locally path connected and, moreover, it is (locally)
a union of paths of increasing distance from $q\in M$.  So if $X$ is such a
neighborhood in a manifold, then
$X \subset \SADH_q(X)$, and $\SADH_q(X)$ is (again) a union of paths 
of increasing distance from $q$.

\begin{remark}
It is easy to see that $\SADH_q(X)$
is contained in the rational hull of $X$:
Consider any rational function $f$ that is holomorphic in a
neighborhood of $X$.
By the Kontinuit\"atssatz,
$f$ analytically continues to any point of $\SADH_q(X)$, and hence
$\SADH_q(X)$ does not intersect the pole set.
\end{remark}

It is useful to know whether the hulls we define are compact if we are starting
with a compact set.  The set $\SADH_q(X)$ is compact if $X$ is, but
unfortunately $\DH(X)$ need not be compact.  However, the polynomial hull is
always compact.

\begin{lemma}
If $X \subset \C^n$ is compact, then $\SADH_q(X)$ is compact.
\end{lemma}

\begin{proof}
That $\SADH_q(X)$ is bounded follows from the maximum principle.
To see that it is closed, suppose $p_1 \in \overline{\SADH_q(X)}$.
For any $\epsilon > 0$ there is a $p_2 \in \SADH_q(X)$ that is within
$\frac{\epsilon}{2}$ of $p_1$.
There is an
approximately attached shrinking family of discs where the last disc goes
through $p_2$ and the boundaries are within $\frac{\epsilon}{2}$ of $X$.
By adding a small function linear in $t$ to the family,
we can create a new family of discs where the last disc goes through $p_1$ and the
boundaries are within $\epsilon$ of $X$.
As $\epsilon > 0$ was arbitrary,
$p_1 \in \SADH_q(X)$.
\end{proof}

\begin{example} \label{example:DHnotcompact}
Given a compact set $X$, the set $\DH(X)$ need not be compact.
Let $\varphi_{k} \colon \overline{\D} \to \C$, $k \in \N$,
be a uniformly bounded sequence of analytic functions such that
$\varphi_k(0)=0$ and 
$\varphi_{k}|_{\partial \D}$ converges to a bounded function that is
continuous except at the point $1 \in \partial \D$.  Assume also
that the sequence $\varphi_k$ converges on $\overline{\D} \setminus \{ 1 \}$ 
to a continuous function $\varphi_{\infty}$ holomorphic on $\D$.
For simplicity, we could also arrange that each $\varphi_k$ extends a little
past the circle.
It is not difficult to construct such a sequence.
Define
\begin{equation}
X = \overline{
\{ z \in \C^3 : ~ z_1 \in \partial \D,
~z_2 = \varphi_k(z_1), ~ z_3 = 1/k, ~ k \in \N \}} .
\end{equation}
Any disc attached to $X$ lies in a set where $z_3$ is constant, and
hence either $z_3 = 0$ or $z_3 = \frac{1}{k}$ for some $k$. For each $k\in\N$ define $\Phi_k(\zeta)=(\zeta, \varphi_k(\zeta), \frac{1}{k})$.
The discs $\Phi_k$ are all attached to $X$, and their images lie in $\DH(X)$.
Moreover, for $k \in \N$, let
$\Psi(\zeta) = \bigl(\alpha(\zeta),\beta(\zeta),\frac{1}{k}\bigr)$ be an analytic disc.
The image of $\Phi_k$ is a subvariety, and it is given by 
functions that are holomorphic a little bit past the boundary. Thus,
the analytic disc $\Psi$ takes the circle into this subvariety,
and hence it takes $\D$ into the subvariety.  In other words, the image 
$\Phi_k(\overline{\D})$ gives all the points in $\DH(X)$ for $z_3 =
\frac{1}{k}$.  Note that $(0,0,\frac{1}{k}) \in \DH(X)$.

Now suppose that
$\Psi(\zeta) = \bigl(\alpha(\zeta),\beta(\zeta),0\bigr)$ is an analytic
disc.
As $\sabs{\alpha(\zeta)} = 1$ for $\sabs{\zeta}=1$, we have that
$\alpha$ is a finite Blaschke product.  In particular, either it is constant,
or it goes around the entire circle some number of times.
If $\alpha$ is constant, then clearly $(0,0,0)$ is not in the image of
$\Psi$.  If $\alpha$ is not constant, then for each $\zeta \in \partial \D$
such that $\alpha(\zeta) \not= 1$, we have
$\beta(\zeta) = \varphi_\infty\bigl(\alpha(\zeta)\bigr)$: the disc is
attached, and that is the only possible $z_2$-coordinate.  In particular,
we find a contradiction as this means that $\beta$ is discontinuous.
Therefore, there are no discs attached to $X$ where $z_3 = 0$
and $z_1$ is not a unimodular constant.  In other words,
$(0,0,0) \notin \DH(X)$.  But $(0,0,0)$ is in the closure
of $\DH(X)$, so $\DH(X)$ is not closed.
\end{example}

\begin{example}\label{example:pert}
The iterated SADH condition is not stable under perturbation.  For example,
consider $M \subset \C^3$  given by $\Im w = \sabs{z_1}^4 - \sabs{z_2}^4$.
The standard technique of attaching discs one normally uses
for a hypersurface with indefinite Levi form applies. These discs fill a
neighborhood of the origin, and all shrink to the origin.
(Affine linear discs suffice.)
However, the perturbation
$M_\epsilon$ given by $\Im w = \epsilon ( \sabs{z_1}^2 + \sabs{z_2}^2) + 
\sabs{z_1}^4 - \sabs{z_2}^4$ is strictly pseudoconvex at $0$ for all $\epsilon > 0$,
and hence all analytic discs attached to $M_\epsilon$ near $0$ must fall on one side of $M_\epsilon$.
In a similar manner, examples having higher codimension can be constructed.
\end{example}

For smooth generic CR submanifolds, local attached analytic discs form a Banach manifold, and all such discs
will shrink to a point; see Sections 6.5 and 8.2 of \cite{BER:book}.  However,
shrinking families of discs are not guaranteed for CR singular submanifolds near the
CR singular point, as the next two examples show.

\begin{example}\label{example:discrete}
Consider the smooth submanifold $M \subset \C^3$ given in coordinates $(z,w_1,w_2)$ via
\begin{equation}
w_1 = \sabs{z}^2 , \quad
w_2 = \sabs{z}^2 + f\bigl(\sabs{z}^2\bigr) (\Re z) ,
\end{equation}
where $f(t)$ is a smooth real-valued function that is zero precisely when
$t=\frac{1}{n}$ for $n\in\N$ or $t=0$.  The submanifold $M$ is of dimension 2. It has a
CR singularity at the origin but is totally real at other points.
Let $\varphi(\zeta) = \bigl(z(\zeta),w_1(\zeta),w_2(\zeta)\bigr)$ be an
analytic disc
attached to $M$. Then $w_1(\zeta)$ and $w_2(\zeta)$ are holomorphic functions that are real-valued  on $\partial \D$, so they are constant on $\overline{\D}$.
This means that $\sabs{z(\zeta)}$ is also constant on $\partial \D$.
If we insist that $\varphi$ be nonconstant, then $z$ itself must be
nonconstant; in particular, $\Re z$ must be nonconstant on $\partial \D$.
But since $\sabs{z}^2 + f\bigl(\sabs{z}^2\bigr) (\Re z)$ must be constant on $\partial \D$,
we have that $f\bigl(\sabs{z}^2\bigr)$ must be zero, which is true only if
$w_1 = \sabs{z}^2 = \frac{1}{n}$.  In other words, the only nonconstant attached
analytic discs to $M$ are those in the discrete sequence of discs
\begin{equation}
\varphi_{n}(\zeta) = \left(\frac{1}{\sqrt{n}} \zeta, \frac{1}{n},
\frac{1}{n} \right).
\end{equation}
This sequence does ``shrink to zero'' discretely but not continuously, so it
does not give a shrinking disc hull, although these discs are in the regular disc hull $\DH(M)$.
In particular, these discs cannot be used via the Kontinuit\"atssatz to extend CR${}_H$ functions beyond the initial neighborhood in which they are defined.
\end{example}

\begin{example}
If we modify the preceding example by taking $f$ to be a real-analytic function with finitely many zeros,
we find a real-analytic submanifold with trivial topology (topology of a ball)
that has only finitely many attached discs.
\end{example}

\begin{example}\label{example:neediterated}
Let us show that iteration may be necessary.
Consider the set $X \subset \C^2$ given by
\begin{multline}
X = X_1 \cup X_2 =
\{ z \in \C^2 : \sabs{z_1} = \sabs{z_2} = 1 \text{ and } \Im z_2 \geq 0 \}
\\
\cup
\{ z \in \C^2 : \sabs{z_1} = 2, \sabs{z_2} = 1, \text{ and } \Im z_2 \leq 0 \} .
\end{multline}
If $\varphi = (\varphi_1,\varphi_2) \colon \overline{\D} \to \C^2$ is an analytic
disc attached to $X$, then as the two components of $X$ are disconnected, we have
that either $\varphi(\partial \D) \subset X_1$ or 
$\varphi(\partial \D) \subset X_2$.
Suppose $\varphi(\partial \D) \subset X_1$.  Then we find that $\varphi_2$
must be constant.  A similar argument applies if $\varphi(\partial \D) \subset X_2$.  In either case, the disc then fills in
all of $\sabs{z_1} \leq 1$ or $\sabs{z_1} \leq 2$.  That is, we find that
\begin{multline}
\DH(X) = 
\{ z \in \C^2 : \sabs{z_1} \leq 1,  \sabs{z_2} = 1, \text{ and } \Im z_2 \geq 0 \}
\\
\cup
\{ z \in \C^2 : \sabs{z_1} \leq 2, \sabs{z_2} = 1, \text{ and } \Im z_2 \leq 0 \} .
\end{multline}
In particular, $\DH(X)$ contains the torus $\bT^2$ given by
$\sabs{z_1}=\sabs{z_2}=1$. However, $\DH(X)$ does not contain the
polydisc $\overline{\D^2}$. It is a relatively routine computation that
$\DH(\bT^2) = \overline{\D^2}$.  Hence $\DH(X)$ does not contain the unit
polydisc, but $\DH^2(X)$ does.
Thus, iteration is necessary for some sets.

The reader may complain that $X$ is disconnected
and that the discs do not all shrink to a point.
We modify the previous example as follows.
Consider the set $X' \subset \C^3$ given by
\begin{multline}
X' = 
\{ z \in \C^3 : \sabs{z_1} = \sabs{z_2} = \Re z_3 , \Im z_2 \geq 0 , z_3 \in[0,1]   \}
\\
\cup
\{ z \in \C^3 : \sabs{z_1} = 2 \Re z_3, \sabs{z_2} = \Re z_3 , \Im z_2 \leq 0, z_3 \in[0,1] \} .
\end{multline}
By $z_3 \in [0,1]$ we mean that $z_3$ is real and in the unit interval.
The set $X'$ is connected and compact. Since the third component of every
analytic disc attached to $X'$ (and hence to any disc hull) must be constant, we
reduce the computation to a scaled version of the above example in $\C^2$.
It is not hard to see that similar reasoning holds also for approximately
attached discs, and moreover
we find that every disc shrinks to the origin.
Hence, we have a
set $X'$ for which $\SADH_0(X')=\DH(X')$ and
$\SADH_0\bigl(\SADH_0(X')\bigr)=\DH^2(X')$. In particular,
if we wish to extend functions via the hull $\SADH_0$, we must iterate:
$\SADH_0\bigl(\SADH_0(X')\bigr) \not= \SADH_0(X')$.
\end{example}

\begin{example}\label{example:nonuiniqueonSADH}
Let us show that a function in $\CR_H(X)$ does not necessarily uniquely
extend to a holomorphic function in a neighborhood of $\SADH_q(X)$ even
though
it admits analytic continuation to such a neighborhood.
Define $X \subset \C^3$ by
\begin{multline}
X =
\left\{ z \in \C^3 :
z_1 = t\zeta,
z_2 = 0, 
z_3 = e^{i\pi t},
\zeta \in \partial \D,
t \in [0,1]
\right\}
\\
\cup
\left\{ z \in \C^3 :
z_1 = 0, 
z_2 = t\zeta,
z_3 = e^{-i\pi t},
\zeta \in \partial \D,
t \in [0,1]
\right\} .
\end{multline}
It is clear that the discs
$\zeta \in \overline{\D} \mapsto 
\left(
t\zeta,0,e^{i\pi t}
\right)$ and
$\zeta \in \overline{\D} \mapsto 
\left(
0,t\zeta,e^{-i \pi t}
\right)$ are attached to $X$, and furthermore these give families 
shrinking to the point $q=(0,0,1)$.  It is therefore clear that
$\SADH_q(X)$ contains the set
\begin{multline}
\left\{ z \in \C^3 :
z_1 = t\zeta,
z_2 = 0, 
z_3 = e^{i \pi t},
\zeta \in \overline{\D},
t \in [0,1]
\right\}
\\
\cup
\left\{ z \in \C^3 :
z_1 = 0, 
z_2 = t\zeta,
z_3 = e^{-i \pi t},
\zeta \in \overline{\D},
t \in [0,1]
\right\} ,
\end{multline}
and any function in $\CR_H(X)$ can be analytically continued along the paths
given by
$t \in [0,1] \mapsto 
\left(
0,0,e^{i \pi t}
\right)$ and
$t \in [0,1] \mapsto 
\left(
0,0,e^{-i \pi t}
\right)$ via the
Kontinuit\"atssatz.  The point $(0,0,-1)$ is in $\SADH_q(X)$ but not in
$X$.  Take a branch of, say, $\sqrt{z_3}$ on some neighborhood of $X$
(which we can do since $X$ is simply connected, and so is some thickening of $X$
to obtain a neighborhood).  This function is in $\CR_H(X)$, but it cannot
possibly extend to a $\CR_H$ function on $\SADH_q(X)$ as the continuation
along the two paths given above will give a different value at $(0,0,-1)$.
\end{example}


\section{Results on hulls and approximations} \label{section:results}

Having an iterated disc hull be a neighborhood at $p$ implies that CR${}_P$
functions of any regularity extend to some neighborhood of
$p$.  

\begin{thm}\label{thm:prophullk}
Let $M \subset \C^n$ be a real submanifold of regularity $C^\ell$
for $\ell \geq 1$ that satisfies the DH condition at $q \in M$.
Then, for $k\leq \ell$, $M$ has the extension property for $\CR_P^k(M)$ at~$q$.
\end{thm}

\begin{proof}
Suppose $f \in \CR_P^k(M)$.  Then there exist a compact neighborhood
$K$ of $q$ in $M$ and a sequence $\{ P_j \}$ of holomorphic polynomials
converging uniformly to $f$ on $K$.
By the maximum principle, $\{ P_j \}$ converges uniformly on $\DH(K)$,
and therefore (by iterating) also on $\DH^N(K)$ for any $N$.  By hypothesis, there exists $N$ such that
$\DH^N(K)$ is a neighborhood of $q$ in $\C^n$, and then
$\{ P_j \}$ converges uniformly on its interior to a holomorphic function extending~$f$.
\end{proof}

Real-analytic CR functions automatically extend to a neighborhood at CR points of a real-analytic submanifold via Severi's theorem.
Thus, we immediately get the following corollary.

\begin{cor}\label{thm:prophullRAk}
Let $M \subset \C^n$ be a real-analytic submanifold 
that satisfies the DH condition at all CR singular points.
Then $\CR_P^\omega(M) = \CR_H(M)$.
\end{cor}

We now want to consider extending $\CR_H$ functions using $\SADH_q$.
Example~\ref{example:nonuiniqueonSADH} says that we cannot just assume that we can
extend functions from $\CR_H(X)$ to $\CR_H\bigl(\SADH_q(X)\bigr)$.  A simple
scenario where we avoid the multi-valuedness issues is when $\SADH_q(X)$ is
already a neighborhood.  This result was proved in the erratum for
\cite{LNR:Cartan}.
For completeness and the reader's convenience,
we reproduce the theorem and the proof here.

\begin{thm}\label{thm:sadhext}
Let $X \subset \C^n$ be a
compact and connected subset and $q \in X$.
Suppose that $\SADH_q(X)$ has $q$ in its interior,
and suppose $B \subset \SADH_q(X)$ is a ball centered at $q$
such that $B \cap X$ is connected.
If $f$ is a holomorphic function defined
on some neighborhood of $X$, then there exists a holomorphic function
$F \colon B \to \C$ such that $f=F$ on $B \cap X$.
\end{thm}

\begin{proof}
Suppose $B = B_{\delta}(q)$.
The function $f$ extends to some $B_{\delta'}(q)$ for $\delta' > 0$
as $f$ is defined in a neighborhood of $X$.
Suppose $\delta'$ is the largest $\delta' \leq \delta$
for which $f$ extends uniquely to
$B_{\delta'}(p)$.  Suppose for a contradiction that $\delta' < \delta$.
We will show that $f$ extends uniquely to a slightly larger ball.
Consider $p \in \partial B_{\delta'}(q)$.
As $p \in \SADH_q(X)$,
there exists a path from $q$ to $p$ of increasing distance from $q$
along which $f$ can be analytically continued. Except for the endpoint $p$, the path lies entirely in $B_{\delta'}(q)$.
There exists a small ball $\tilde{B}$ centered at $p$ such that
$f$ extends uniquely to $B_{\delta'}(q) \cup \tilde{B}$.
This construction can be done at every point in
$\partial B_{\delta'}(q)$, which is 
compact.
Thus, $f$ extends uniquely to some slightly larger $B_{\delta''}(q)$.
\end{proof}

Note that every submanifold is locally simply connected.
If, in the context of Theorem \ref{thm:sadhext}, $X=M$ is a submanifold that is simply connected or generic
at CR points, then every $\CR_H(M)$ function extends to $B$.
Because for a submanifold we can
suppose that $B$ is small enough so that $B \cap M$ is connected,
we have the following corollary.

\begin{cor} \label{cor:bigSADHnoiter}
Let $M \subset \C^n$ be a $C^1$ real submanifold that is simply connected or
generic at CR points. If $q \in M$
and
$\SADH_q(M)$ contains a neighborhood (in $\C^n$) of $q$,
then $M$ has the fixed-neighborhood extension property for $\CR_H(M)$ at $q$.
\end{cor}

As we saw, sometimes iteration is necessary to obtain a neighborhood.
In this case, we will need the more complicated iterated SADH condition to get
the fixed-neighborhood extension property for $\CR_H$.

\begin{thm} \label{thm:prophullSADH}
If $M \subset \C^n$ is a $C^1$ real submanifold that is
 simply connected or generic at all CR points and
satisfies the iterated SADH condition at $q \in M$,
then $M$ has the fixed-neighborhood extension property 
for $\CR_H(M)$ at $q$.
\end{thm}

\begin{proof}
Let $f \in \CR_H(M)$.  By the assumptions on $M$, $f$ extends uniquely to some
neighborhood of $M$ as a holomorphic function.
Supposing that $M$ has SADH iterates $X_0 = M$, $X_1$, \ldots, $X_k$ at $q$, the 
proof is then an iteration of Lemma~\ref{lemma:prophullSADH},
as the lemma applies to each $X_j$ to get an extension to a neighborhood of
$X_{j+1}$.  Because $X_k$ contains a neighborhood of the origin (which depends
only on $M$ and not on $f$), $M$ has the fixed-neighborhood extension
(and hence approximation) property for $\CR_H(M)$ at $q$.
\end{proof}

We now show that, given a set $X$, if we take $Y$ to be
the union of a compact family of nonintersecting paths from $\SADH_q(X)$ along which we have
analytic continuation as in the definition of
the iterated SADH condition, then we have a uniquely defined holomorphic function on a neighborhood of
$\SADH_q(X)$.

%
%

\begin{lemma} \label{lemma:prophullSADH}
Suppose $q \in \C^n$ and $\sY$ is a compact (with respect to the uniform norm)
set of paths $\psi \colon [0,1] \to \C^n$
with $\psi(0)=q$ and $\snorm{\psi(t)-q}$
strictly increasing in $t$,
and whenever
$\psi_1,\psi_2 \in \sY$ are such that
$\psi_1(t_1)=\psi_2(t_2)$ for some $t_1,t_2 > 0$,
then $\psi_1([0,t_1]) = \psi_2([0,t_2])$.
Suppose further that $f$ is a holomorphic function defined near $q$
that may be analytically continued along all $\psi \in \sY$.  Then
there exists a holomorphic function $F$ defined in a neighborhood of $Y$,
where
$Y = \bigcup_{\psi \in \sY} \psi\bigl([0,1]\bigr)$ is the image of
$\sY$,
such that $f$ and $F$ agree near $q$.
\end{lemma}

\begin{proof}
As $f$ is defined in some neighborhood of $q$, there must be
some positive $\delta > 0$ such that a function $F$ well-defined in a
neighborhood of $Y \cap B_\delta(q)$ exists and agrees with $f$ near $q$.  The lemma will follow if we show that $F$ can be extended into a larger
ball.  The function $f$, and hence $F$, can be analytically continued along the paths in
$\sY$ to a larger set by assumption, so what we need to show is
single-valuedness.

For every point $p \in Y \cap \partial B_\delta(q)$ the function
$F$ can be continued through $p$ along a path of strictly increasing
distance, and hence there is a small ball $B_\epsilon(p)$
with a continuation defined on $B_\epsilon(p)$ and agreeing with
$F$ on a neighborhood of $B_\epsilon(p) \cap Y \cap B_\delta(q)$.
Suppose for a contradiction that for every small enough $\epsilon > 0$ there
are two points $p_1$ and $p_2$ in $Y \cap \partial B_\delta(q)$
such that $B_\epsilon(p_1) \cap B_\epsilon(p_2) \not= \emptyset$
and the continuations in $B_\epsilon(p_1)$ and
$B_\epsilon(p_2)$ do not agree on the intersection.
Each such point corresponds to a path, and thus there must exist
a sequence of pairs of such paths which (after possibly
passing to a subsequence) converge in $\sY$ to a fixed path $\psi$ through
some $p$.
There must therefore be a continuation along $\psi$ to an $\epsilon$-ball
around $p$ as above, but this would be a contradiction to the existence
of the above sequence.

Thus, there is a single $\epsilon > 0$, such that $F$ continues
as above to $B_\epsilon(p)$, and whenever two
such balls intersect, the continuation agrees on the intersection.
Hence $F$ has a single-valued extension to some neighborhood of
$Y \cap B_{\delta'}(q)$ for some $\delta' > \delta$ and we are
done.
\end{proof}

These ideas may be combined in various ways, and at this point we give one
such corollary of the proof.
We will find it useful to extend the proof to submanifolds
$M \subset \C^n \times \R$. 
In this case any SADH iterate is a subset of
$\C^{n} \times \R$, and therefore we will not find
a neighborhood in $\C^{n+1}$.
However, we may find
a neighborhood in the topology of $\C^n \times \R$
to which we can extend all CR${}_H$ functions via the lemma.
Here is the formal statement.

\begin{cor}\label{cor:prophullSADHflat}
Let $M \subset \C^n \times \R$ be a $C^1$ real submanifold
that is simply connected or generic at all CR points,
and let $q \in M$.
Suppose $X_0=M$ and there exist SADH iterates
$X_1,\ldots,X_k$,
where $X_k$
contains a neighborhood of $q$ in $\C^n \times \R$.
Then there exists a neighborhood $U \subset \C^n \times \R$
of $q$
such that every function in $\CR_H(M)$ extends to
a function in $\CR_H(U)$.
\end{cor}

If $M$ is real-analytic, the CR${}_P$ functions that are real-analytic extend at all CR points by Severi's theorem. If $M$ also satisfies the strong iterated SADH condition at the CR singular points, it turns out that we get the fixed-neighborhood extension and approximation properties
for $\CR_P^\omega$ at such points.  We first need a lemma to show that the points that can be achieved by approximately attached discs are
also in the polynomial hull.

\begin{lemma} \label{lemma:approxdiscinPhull}
Suppose that $X \subset \C^n$ is a compact set and $p \in \C^n$
is such that, for every $\epsilon > 0$, there exists an approximately
attached analytic disc $\varphi \colon \overline{\D} \to \C^n$
such that $\varphi(\partial \D) \subset X_\epsilon$ and $p\in \varphi(\D)$.
Then for every holomorphic polynomial $P \colon \C^n \to \C$ we have
\begin{equation}
\sabs{P(p)} \leq \sup_{z \in X} \sabs{P(z)} .
\end{equation}
\end{lemma}

\begin{proof}
Fix $P$ and $\delta > 0$.  As $X$ is compact, there exists $\epsilon > 0$
such that 
\begin{equation}
\sup_{z \in X_\epsilon} \sabs{P(z)} \leq \sup_{z \in X} \sabs{P(z)} + \delta .
\end{equation}
Now pick a $\varphi$ that is $\epsilon$-approximately attached to $X$ and use
the maximum principle:
$\sabs{P(p)} \leq \sup_{z \in X} \sabs{P(z)} + \delta$.
\end{proof}


We note that the strong iterated SADH condition at a point implies the DH condition,
and we get the following corollary.  We require the strong
SADH condition as we want discs that are attached arbitrarily near the point
$q$.  Although these discs all shrink, we do not have control on how fast they
shrink for different paths.

\begin{cor}\label{thm:prophull}
Let $M \subset \C^n$ be a real-analytic submanifold such that
$M$ satisfies the strong iterated SADH condition at every CR singular point of $M$.
Then if $q \in M$ is a CR singular point,
$M$ has the fixed-neighborhood extension property for $\CR_P^\omega(M)$ at $q$,
and thus the fixed-neighborhood approximation property for $\CR_P^\omega(M)$ at $q$.
\end{cor}

\begin{proof} First we show that $\CR_P^\omega(M)\subset \CR_H(M)$.
Suppose $f \in \CR_P^\omega(M)$.   At each CR point $p \in M$,
$f$ extends to a neighborhood of $p$ as a holomorphic function
because $f$ is real-analytic. Now assume that we have the strong iterated SADH condition at $p\in M$.
By Lemma \ref{lemma:approxdiscinPhull} the approximately attached discs are still
in the polynomial hull. Thus, 
given a sequence of polynomials converging uniformly on some neighborhood $X_0$
of $p$ in $M$, we get convergence of this sequence on all the sets
$X_j$ from Definition~\ref{defn:SADHcondition}: Use the same argument
as in the proof of Theorem~\ref{thm:prophullk}, but now use
Lemma~\ref{lemma:approxdiscinPhull}. Therefore, the sequence
converges on a
neighborhood of $p$ in $\C^n$.  Thus, $f$ extends to a neighborhood of $p$. Hence, $f$ extends to a neighborhood of each CR singular point of $M$. We conclude that $f \in \CR_H(M)$.

Now apply Theorem~\ref{thm:prophullSADH} on a simply connected
neighborhood of $q$ to obtain the result.
\end{proof}

If one can extend CR functions in some way near CR singular
points, then one can use other techniques
to extend at the CR points.
In the following corollary we use the Hanges--Tr\`eves
theorem
to propagate
the extension property along complex curves.

\begin{cor}\label{thm:HangesTreves}
Let $M \subset \C^n$ be a smooth real submanifold such that
through every CR point of $M$ there is a connected nonsingular complex curve
$C \subset M_{CR}$ such that the closure of $C$ contains a CR singular point of $M$.
Suppose $f \in \CR^k(M)$ for some $k$ and that,
near every CR singular point, $f$ is the restriction of a holomorphic
function on a neighborhood in $\C^n$.  Then $f \in \CR_H(M)$.

In particular, if in addition $M$ satisfies the DH condition
at each CR singular point, then $\CR_P^k(M) = \CR_H(M)$  for all $k$.

Moreover, if in addition at some $p \in M$, a simply connected
neighborhood of $p$ in $M$ satisfies the iterated SADH condition at $p$, then
$M$ has the fixed-neighborhood extension property
(and thus the fixed-neighborhood approximation property) for $\CR_P^k(M)$
at $p$.
\end{cor}

\begin{proof}
Suppose $f \in \CR^k(M)$, $q$ is a CR singular point, $f$
extends holomorphically to a neighborhood of $q$, and $C \subset M_{CR}$
is a connected nonsingular complex curve
 whose closure contains $q$. Then $f$
extends holomorphically to a neighborhood of at least one
point of $C$.
By the theorem of Hanges and Tr\`eves (Theorem 4.1 of \cite{HangesTreves}),
$f$ extends holomorphically to a neighborhood at each point of $C$.
The first part follows.

The rest of the corollary follows by applying
Theorems \ref{thm:prophullk} and \ref{thm:prophullSADH}.
\end{proof}

As we noted earlier, there are other ways to combine these ideas, but it doesn't
seem productive to list all of the different possibilities;
we have listed only those that seem most relevant for our purposes.

Constructing disc hulls is generally easier
than constructing shrinking disk hulls.
Homogeneity of the set allows us to pass from
disc hulls to shrinking disc hulls.
We say that a set $X \subset \C^n$ is \emph{bounded weighted homogeneous}
if $X$ is bounded and
there exists $\alpha \in \N^n$
such that if $z \in X$ then
$(t^{\alpha_1} z_1, \ldots, t^{\alpha_n} z_n) \in X$ for all $t \in [0,1]$.
For bounded weighted homogeneous sets,
we can always consider all of $X$ rather than an arbitrary
neighborhood of the origin, as we can always rescale our variables.

\begin{lemma}\label{lemma:rein}
Suppose $X \subset \C^n$ is a bounded weighted homogeneous compact set.
Then the sets $X_0=X$, $X_1 = \overline{\DH^1(X)}$, \ldots, $X_{\ell} =
\overline{\DH^{\ell}(X)}$ are SADH iterates at $0$.
If there exists $k$ such that $\overline{\DH^k(X)}$
contains a nonempty Reinhardt domain $V$,
then $X_{k+1} = \overline{\DH^{k+1}(X)}$
contains a complete Reinhardt domain containing $V$.
Therefore, $X$ satisfies the strong iterated SADH condition at $0$.
\end{lemma}

\begin{proof}
By definition, there exists $\alpha \in \N^n$ such that, with $\delta_t$
defined by $\delta_t(z)=(t^{\alpha_1} z_1, \ldots, t^{\alpha_n} z_n)$, we
have $\delta_t(X)\subset X$ for $t\in [0,1]$.
The same $\alpha$ applies to all bounded weighted homogeneous sets in the
argument that follows. 

{\it Claim:} $\overline{\DH(X)}$ is bounded weighted homogeneous,
and $\overline{\DH(X)} \subset \SADH_0(X)$.
Moreover, $\overline{\DH(X)}$ is a 
union of a compact family of nonintersecting paths from $\SADH_0(X)$
as in Definition~\ref{defn:SADHcondition},
that is, it is a SADH iterate at $0$.

{\it Proof:}
By the maximum principle, $\DH(X)$ is bounded because $X$ is bounded. 
Fix $z\in \DH(X)$. We prove that $\delta_t(z)\in \DH(X)$ when $t\in[0,1]$
and that $z\in\SADH_0(X)$.
By definition, there exist $\zeta\in\D$ and an
analytic disc $\varphi \colon \overline{\D} \to \C^n$ such that
$\varphi(\zeta) = z$ and $\varphi(\partial \D) \subset X$. 
Fix $t\in [0,1]$ and define $\varphi_t=\delta_t\circ \varphi$. 
Then $\varphi_t$ is an analytic disc, $\varphi_t(\zeta)=\delta_t(z)$, and 
$\varphi_t(\partial \D) \subset \delta_t(X)$; thus, $\delta_t(z)\in \DH(\delta_t(X))\subset \DH(X)$.
Further, $\varphi_1(\zeta)=\varphi(\zeta) = z$
and $\varphi_0 \equiv 0$. Thus, $z\in\SADH_0(X)$.
For later use, we remark that we showed $\delta_t(\DH(X))\subset\DH(\delta_t(X))$.

We have proved that $\DH(X)$ is bounded weighted homogeneous and 
that $\DH(X)\subset\SADH_0(X)$.  It follows that $\overline{\DH(X)}$ is also
bounded weighted homogeneous, and
$\overline{\DH(X)}\subset\SADH_0(X)$ since $\SADH_0(X)$ is compact.

The mapping $\delta_t$ gives nonintersecting paths in $\overline{\DH(X)}$
which come from $\SADH_0(X)$ as required.
Because $\overline{\DH(X)}$ is compact, the set of paths given by
$\delta_t(z)$ for $z \in \overline{\DH(X)}$ is compact in the uniform norm.
The claim is proved.

By iterating the claim we find that each $\overline{\DH^\ell(X)}$ is bounded weighted
homogeneous and a union of a compact family of nonintersecting paths from
$\SADH_0\bigl(\overline{\DH^{\ell-1}(X)}\bigr)$.
Therefore, these sets are candidates for use in
Definition~\ref{defn:SADHcondition}.  What is missing is to show that
$\overline{\DH^{k+1}(X)}$ contains a complete Reinhardt domain if $\overline{\DH^k(X)}$
contains a nonempty Reinhardt domain.


{\it Claim:} Assume that $X$
contains a nonempty Reinhardt domain $V$.
Then $\overline{\DH(X)}$ 
contains a complete Reinhardt domain containing $V$.

{\it Proof:}
Suppose $X$ contains a nonempty Reinhardt domain $V$.
Define $W=\cup_{t\in [0,1]}\delta_t(V)$.
Then $W$ is bounded weighted homogeneous and invariant under
the action of the $n$-dimensional torus, and $W\subset \DH(X)$. 
Pick a polyradius $(r_1,\ldots,r_n)$ such that, if $\sabs{z_j} = r_j$
for all $j$, then $z \in W$.  All of the discs
\begin{multline}
r_1 \overline{\D} \times r_2 \partial \D \times \cdots \times r_n \partial \D,\quad
r_1 \partial \D \times r_2 \overline{\D} \times r_3 \partial \D \times \cdots \times r_n \partial \D,\quad
\ldots,
\\
r_1 \partial \D \times r_2 \partial \D \times \cdots
\times
r_{n-1} \partial \D \times r_n \overline{\D},
\end{multline}
are attached to $W$.
As $W$ is bounded weighted homogeneous, we find that these discs
composed with $\delta_t$ for
$t\in[0,1]$ are also attached to $W$, and hence all of these
discs are in $\DH(X)$.
Therefore, the entire polydisc of polyradius $(r_1,\ldots,r_n)$
is in $\DH(X)$.  The claim is proved.

By iterating the claim, as each $\overline{\DH^\ell(X)}$ is bounded weighted homogeneous
if $X$ is, we obtain the conclusion of the lemma about $X_{k+1}$. 

It remains to
show that $X$ satisfies the strong iterated SADH condition at $0$.  This fact
follows by bounded weighted homogeneity.  For every neighborhood $U$ of the origin,
there is some $t \in (0,1]$ such that $\delta_t(X) \subset U$.  Moreover,
as every set above is weighted homogeneous, the conclusions of the claims
also follow for $\delta_t(X) \subset X$, and the desired conclusion follows.
We remark that in this setting, the number $k$ is independent of the neighborhood
$U$.
%
%
%
\end{proof}

We noted earlier that, when $M$ is a subset of
$\C^{n} \times \R \subset \C^{n+1}$,
it cannot satisfy the iterated SADH condition at any point.
However, such manifolds can have the
fixed-neighborhood approximation property by the following
generalization of the Weierstrass approximation theorem.
One such example is
$w=|z_1|^2-|z_2|^2$ (see Section~\ref{section:BTnofixednbhd}).


\begin{thm}\label{thm:graph}
Let $M \subset \C^n \times \R \subset \C^{n+1}$ be a (topological) submanifold
given as a graph $s=\rho(z,\bar{z})$, where
$(z,s) \in \C^n \times \R$ denote the variables and $\rho$ is continuous.
Suppose $K \subset \C^n \times \R$ is a compact neighborhood of $q \in M$
and $M \cap K$ is a nonempty compact set.
Let $K_s = \{ z \in \C^n : (z,s) \in M \cap K \}$, and suppose that
\begin{enumerate}[(*)]
\item
For any $\epsilon > 0$ there is a $\delta > 0$ such that if
$\sabs{s-t} < \delta$ and $K_s$ and $K_t$ are nonempty,
then $d_H(K_s,K_t) < \epsilon$, where $d_H$ denotes the Hausdorff
distance.
\end{enumerate}
Let $\sF$ be a class of continuous functions on $M$ with the following property:
For every $f \in \sF$ and
each $s \in \R$ for which
$K_s$ is nonempty,
the function $K_s\ni z \mapsto f(z,s)$ can be uniformly
approximated on $K_s$ by polynomials in $z$.

\nopagebreak
Then every $f \in \sF$ can be uniformly approximated on $M \cap K$ by
polynomials in $(z,s)$.
\end{thm}

\begin{remark} \label{rmk:graph}
If $\rho \in C^3$, $\nabla \rho|_q = 0$,
and the Hessian of $\rho$ at $q$ is
nondegenerate, then (*) is satisfied for a small enough compact neighborhood $K$.
To see this fact, suppose $q=0$
and apply the Morse lemma
to find a $C^1$ (not holomorphic) change of the $z$ variables near $0$
so that $M$ is given by an equation of the form
$s = \sum_{k=1}^{2n} \pm x_k^2$, where $x_1,\ldots,x_{2n}$ are the new coordinates
for $\C^n$.  The condition (*)
is clearly satisfied in this setting for small enough $K$.  As a $C^1$
diffeomorphism will leave the Hausdorff distance locally comparable, the condition (*) is also true before the change of variables for small enough $K$.
\end{remark}

\begin{proof}[Proof of Theorem \ref{thm:graph}]
Let $f \in \sF$ and $\epsilon > 0$ be given.
Given $s_0$ for which $K_{s_0}$ is nonempty, find
a holomorphic polynomial $P_{s_0}(z)$ such that $P_{s_0}$ is within $\epsilon$ of
$z \mapsto f(z,{s_0})$ on $K_{s_0}$.
Then there exists a neighborhood of $K_{s_0}$ on which $P_{s_0}(z)$ is
within $3\epsilon$ of $f\bigl(z,\rho(z)\bigr)$ because both $P_{s_0}$ and
$f\bigl(z,\rho(z)\bigr)$ are uniformly continuous on a neighborhood of
$\pi_1(M \cap K)$
(the projection of $M \cap K$ onto the $z$-coordinate).
Via the hypothesis (*), there exists $\delta_{s_0} > 0$ such that
$P_{s_0}(z)$ and $f\bigl(z,\rho(z)\bigr)$ are within $3\epsilon$ when $\sabs{s-s_0} < \delta_{s_0}$ and $z\in K_s$.

The set $I = \{ s \in \R : K_s \neq \emptyset \}$ is compact.  So there
exist $s_1 < \cdots < s_\ell$ such that
$s_j \in I$ for each $j$ and the intervals
$I_j = \bigl(s_j-\delta_{s_j},s_j+\delta_{s_j}\bigr)$ cover
$I$. Let $\{\varphi_j\}$ be a continuous partition of unity on $I$ subordinate to $\{I_j\}$, so $\sum_j \varphi_j=1$ on $I$ and for all $j$ we have $\varphi_j\geq 0$ and ${\supp } \; \varphi_j\subset I_j$.
If $\varphi_j(s)\neq 0$, then $s\in I_j$, and thus $|P_{s_j}(z)-f(z,\rho(z))|<3\epsilon$ for $z\in K_s$. Define $P(z,s)=\sum_j \varphi_j(s) P_{s_j}(z)$. If $s\in I$ and $z\in K_s$, we have
\begin{equation}
|P(z,s)-f(z,\rho(z))| = \biggl|\sum_{\{j\colon \varphi_j(s)\neq 0\}} \varphi_j(s) [P_{s_j}(z)-f(z,\rho(z))]\biggr|<3\epsilon \sum_{\{j\colon \varphi_j(s)\neq 0\}}\varphi_j(s) =3\epsilon.
\end{equation}

Thus, $P\bigl(z,\rho(z)\bigr)$ is within $3 \epsilon$ of $f\bigl(z,\rho(z)\bigr)$ for every $z \in \pi_1(M \cap K)$. Also,
$P(z,s)$ is a polynomial in $z$, and if we write
\begin{equation}
P(z,s) = \sum_{\alpha} a_{\alpha}(s) z^\alpha
\end{equation}
then the coefficients $a_\alpha$ are continuous functions on $I$. These functions can be uniformly approximated on $I$ by polynomials in $s$ using the standard
Weierstrass approximation theorem.  By choosing a sufficiently close
approximation, we find a polynomial $Q(z,s)$ such that
$Q(z,s)$ is within $\epsilon$ of $P(z,s)$ for $s\in I$ and $z\in K_s$.  Then
$Q\bigl(z,\rho(z)\bigr)$ is within $4\epsilon$ of $f\bigl(z,\rho(z)\bigr)$
on $\pi_1(M \cap K)$.
\end{proof}

\begin{example}
The hypothesis on the Hausdorff distance in Theorem \ref{thm:graph} is sufficient but not necessary.
Consider $M \subset \C \times \R$ defined by
\begin{equation}
s =
\begin{cases}
e^{-1/(\Re z)^2} & \text{ if } \Re z > 0, \\
0 & \text{ if } \Re z \leq 0 .
\end{cases}
\end{equation}
Then condition (*) of the theorem is not
satisfied: no matter how small a neighborhood $K$ of the origin we take,
$K_0 = \{ z : \Re z \leq 0 , (z,0) \in K \}$
is not contained in a small neighborhood
of $K_s = \{ z : \exp(-1/(\Re z)^2) = s, (z,s) \in K \}$ for
any $s > 0$.

Suppose for simplicity that $K$ is the unit
polydisc.  Any function on $M$
can be represented as a function of $z$ because
$M$ is a graph.
If $f$ is a continuous function on $M$
that can be uniformly approximated by
holomorphic polynomials in $z$ on $K_s$ for 
every $s$, we know that $f$ must
be holomorphic for $\Re z < 0$.

The proof that $f$ can be approximated by a polynomial in $(z,s)$ on $K \cap M$
follows the same logic as the proof of Theorem~\ref{thm:graph} except we start
with $s_1=0$, and we let $P_0$ be a polynomial approximating $f$ for
$\Re z \leq 0$.  The rest of the proof works as before simply considering
the set of $M \cap K$ where $\Re z \geq 0$.
\end{example}


\section{Flat elliptic Bishop surfaces} \label{section:BishopFlatElliptic}

The results in this section concern flat Bishop surfaces that are elliptic. 
With regard to the fixed-neighborhood approximation property for $\CR_P^\omega(M)$, in the special case $w=|z|^2$ we show in Section \ref{section:BTBishop} how to produce approximants using an integral formula for a larger class of functions.

\begin{thm}\label{thm:elliptic}
\leavevmode
Fix $\lambda\in [0,1/2)$. Define  $\rho(z,\bar{z})=z\bar{z}+\lambda(z^2+\bar{z}^2)+E(z,\bar{z})$, where $E(z,\bar{z})$ is smooth, real-valued, and $O(|z|^3)$.
For $\delta_1>0$ and $\delta_2>0$ sufficiently small, define $M=\{(z,w)\in \C^2\colon w=\rho(z,\bar{z}), |z|<\delta_1, |w|<\delta_2\}$. 

\begin{enumerate}[(i)] 
\item If $M$ is $C^\ell$, for all $k\leq \ell$ we have $\CR^k(M) \supsetneq \CR^k_P(M)$, and hence $\CR^k(M)$ does not have the fixed-neighborhood approximation property at the origin.
\item $\CR_H(M)$ does not have the fixed-neighborhood extension property at the origin.  
\item
If $M$ is real-analytic we have:
\begin{enumerate}[(a)]
\item $\CR^\omega_P(M)=\CR_H(M)$
\item $\CR_P^\omega(M)$ has the fixed-neighborhood approximation property at the origin.
\end{enumerate}

\item
In the smooth category:
\begin{enumerate}[(a)]
\item
$\CR^\infty_P(M)\supsetneq\CR_H(M)$
\item
$\CR^\infty_P(M)$ does not have the fixed-neighborhood approximation property at the origin.  
\end{enumerate}

\end{enumerate}
\end{thm}

\begin{proof}

Put $\ell=\omega$ if $E$ is $C^\omega$ and $\ell=\infty$ otherwise.

By \cite{crext1}*{Proposition 3.1} there exist $\delta_1>0$ and $\delta_2>0$ such that the following hold: For every $s\in (0,\delta_2)$, the set \begin{equation}
    K_s=\{z\in\C\colon |z|<\delta_1, s=\rho(z,\bar{z})\}
\end{equation}
is either empty or a connected compact real curve homeomorphic to a circle, and  $K_s$ bounds a relatively compact domain $\Omega_s$ with connected boundary. Moreover, if we define \begin{equation}
   M=\{(z,w)\in \C^2\colon w=\rho(z,\bar{z}), |z|<\delta_1, |w|<\delta_2\},\end{equation} 
then the origin is the only CR singular point of $M$, and $M$ is totally real away from the origin. It follows that $C^k(M)=\CR^k(M)$ for every $k\leq \ell$.

  By \cite{crext1}*{Lemma 6.1}, a smooth function $f$ on $M$ has a holomorphic extension on each nonempty leaf (i.e., a continuous extension from the set $K_s$ that is holomorphic on $\Omega_s$) if and only if the following moment condition holds:
 for each $t>0$ and $k \in \N\cup \{0\}$, we have \begin{equation} \label{eq:moment}
 \int_{M\cap\{w=t^2\}} f(\zeta)\zeta^k \; d\zeta=0.  \end{equation}
Here we think of $M$ as being parametrized by $z$ and consider the corresponding function $f(z)$.
We remark that, if $f\in\CR^0_P(M)$, then  $f$ satisfies the moment condition on some neighborhood of the origin.

Fix $k\leq \ell$. The function $f(z,w)=\bar{z}$, considered as a function on $M$, belongs to $C^\ell (M)\subset\CR^k(M)$. We claim that, on all neighborhoods of the origin, the moment condition fails to hold for $f$. Suppose that the claim is false. Then, by \cite{crext1}*{Theorem 1.1}, there exist small $\tilde{\delta}_1>0 $ and $\tilde{\delta}_2>0$ such that $f$ can be extended to be a smooth function on   \begin{equation}
    \{(z,w)\in \C^2\colon \Re w \geq \rho(z,\bar{z}), \Im w=0, |z|<\tilde{\delta}_1, |w|<\tilde{\delta}_2\}.\end{equation}
Moreover, the extension has a formal power series in $z$ and $w$ at the origin. This is impossible, so we have a contradiction.  Thus, the claim holds, and it follows from the preceding remark that $f\not\in \CR^k_P(M)$.
Hence, for every $k\leq \ell$, $\CR^k(M) \neq \CR^k_P(M)$, and $\CR^k(M)$ does not have the fixed-neighborhood approximation property at the origin.

Next we assume that $\ell=\omega$ and consider functions that are real-analytic on $M$. Note that $\CR_H(M)$ does not have the fixed-neighborhood extension property at the origin because $M$ is contained in the
Levi-flat hypersurface given by $\Im w = 0$.
Now we show that $\CR^\omega_P(M)=CR_H(M)$. Fix $f\in\CR^\omega_P(M)$. As we remarked earlier, $f$ satisfies the moment condition on a neighborhood of the origin.
By \cite{crext1}*{Theorem 1.1}, $f$ extends to be holomorphic on a neighborhood of the origin.  Because the origin is the only CR singular point of $M$,   $f\in CR_H(M)$. This holds for all $f\in\CR^\omega_P(M)$, so $\CR^\omega_P(M)=CR_H(M)$. 

Next we show that  $\CR^\omega_P(M)$ has the fixed-neighborhood
approximation property at the origin. Fix $\tilde{\delta}_1>0$ and $\tilde{\delta}_2>0$ sufficiently small. Let $f\in\CR^\omega_P(M)$. Because $f$ satisfies the moment condition on a neighborhood of the origin, 
and that condition involves the vanishing of certain functions that are real-analytic on the interval $(0,\delta_2)$,
the moment condition holds on all of $M$. Thus, we can find a holomorphic extension of $f$ from each nonempty $K_s$ to $\Omega_s$.
We then use Mergelyan's theorem to uniformly approximate $K_s\ni z \mapsto f(z,s)$ by polynomials in $z$. (Here we use the fact that $\C\setminus\overline{\Omega_s}$ is connected.) This holds whenever $K_s$ is nonempty, so we may use Theorem \ref{thm:graph} (along with Remark~\ref{rmk:graph}) to uniformly approximate $f$ on $M\cap \{(z,s): |z|\leq \tilde{\delta}_1,  |s| \leq \tilde{\delta}_2\}$ by polynomials in $(z,s)$. To get holomorphic polynomials on $\C^2(z,w)$, replace $(z,s)$ by $(z,w)$. This proves that $\CR^\omega_P(M)$ has the fixed-neighborhood
approximation property at the origin. 

Now we consider the class $C^\infty(M)$. We prove that
$\CR^\infty_P(M)$ does not have the fixed-neighborhood
approximation property at the origin. (It then follows that $\CR^\infty_P(M)\neq CR_H(M)$.) Let $\epsilon>0$ be sufficiently small relative to $\delta_1,\delta_2$. Choose $\chi_\epsilon\colon [0,\infty)\rightarrow [0,\infty)$ to be smooth and satisfy $\chi\equiv 0$ on $[0,\epsilon]$ and  $\chi\equiv 1$ on $[2\epsilon,\infty)$. Define $f_\epsilon$ on $M$ by $f_\epsilon(z,\rho(z,\bar{z}))=\chi(|z|^2)\bar{z}$. Then $f_\epsilon\in \CR^\infty_P(M)$ because $f_\epsilon$  is identically $0$ near the origin, and at CR points 
we can apply the Baouendi--Tr\`eves approximation theorem. But for every compact neighborhood of the origin there exists $\epsilon$ such that $f_\epsilon$ does not satisfy the moment condition on that neighborhood. (Use the preceding argument that, on all neighborhoods of the origin, the moment condition fails to hold for $\bar{z}$.) Thus, $\CR^\infty_P(M)$ does not have the fixed-neighborhood
approximation property at the origin.
\end{proof}


\section{Baouendi--Tr\`eves for a special elliptic Bishop surface} \label{section:BTBishop}

Define $M=\{(z,w)\in\C^2\colon w=|z|^2\}$. 
In this section, we prove that the class of functions in $C^0(M)$ satisfying the moment condition, equation \eqref{eq:moment} from Section \ref{section:BishopFlatElliptic}, on a fixed neighborhood of the origin has the fixed-neighborhood
approximation property at the origin. (In fact, for the class $C^\omega(M)$, if the moment condition holds for, say, a nonempty open interval of values of $t$, then it holds for all $t$. See the argument in Section \ref{section:BishopFlatElliptic}.)  The proof produces approximants by means of an integral formula, and in that way it is similar to the original proof of the Baouendi--Tr\`eves approximation theorem.

We make a couple of preliminary comments. Fix $f\in C^0(M)$, and 
for each $t>0$ write the value of $f$ at $(z,t^2)\in M$ as $f_t(z)$.  First, recall from Section  \ref{section:BishopFlatElliptic} the remark that a necessary condition for $f$ to belong to $\CR^0_P(M)$ is that the moment condition hold for $t>0$ sufficiently small.
Second, note that this moment condition is equivalent to
\begin{equation} \label{eq:moment2}
\int_{0}^{2\pi} f_t(te^{i\theta})e^{i(k+1)\theta} \, d\theta=0
\avoidbreak
\end{equation}
for $t>0$ sufficiently small and for all $k\in\N\cup\{0\}$.

Now we consider the class of functions in $C^0(M)$ satisfying the moment condition on a fixed neighborhood of the origin. We prove using an integral formula that this class has the fixed-neighborhood approximation property at the origin.
Let $\epsilon>0$ be given. Choose a nonnegative smooth function $\chi$ on $[0,\infty)$ such that $\chi\equiv 1$ on $[0,\epsilon/2]$ and $\chi\equiv 0$ on $[\epsilon,\infty)$.
 For $n\in \N$ define $c_n$ by $1/c_n=\int_{\C}  \exp{(-|\zeta|^2/n)} \, dA(\zeta)$. 
 
 Given a continuous function on $M$ satisfying the moment condition on $\{(z,|z|^2)\colon |z|\leq \epsilon\}$, we think of $M$ as being parametrized by $z$ and consider the corresponding function $f(z)$.  
 Define 
 \begin{equation} \begin{split} Q_n(z,\bar{z}) & = c_n\int_\C \chi(|\zeta|) f(\zeta) \exp{(-|z-\zeta|^2/n)} \, dA(\zeta)\\
 & = c_n\int_0^\infty\int_0^{2\pi} \chi(r) f(re^{i\theta}) \exp{(-|z-re^{i\theta}|^2/n)}\, r \, d\theta dr.
  \end{split} \end{equation}
 Now  \begin{equation}  \exp{(-|z-\zeta|^2/n)}=\sum_{k=0}^\infty \frac{(-1)^k}{n^kk!}  (z\bar{z}-z\bar{\zeta}-\bar{z}\zeta+\zeta\bar{\zeta})^k, \end{equation} and we write
 \begin{equation} (z\bar{z}-z\bar{\zeta}-\bar{z}\zeta+\zeta\bar{\zeta})^k=\sum a^{(k)}_{\alpha \beta\gamma\delta}(z\bar{z})^\alpha(z\bar{\zeta})^\beta(\bar{z}\zeta)^\gamma(\zeta\bar{\zeta})^\delta. \end{equation} Then 
 $\int_0^{2\pi} f(re^{i\theta}) \exp{(-|z-re^{i\theta}|^2/n)} \, d\theta$ can be written as a sum of constant multiples of terms of the form
 \begin{equation}  \int_0^{2\pi} f(re^{i\theta})(z\bar{z})^\alpha(zre^{-i\theta})^\beta(\bar{z}re^{i\theta})^\gamma r^{2\delta}  \, d\theta= (z\bar{z})^\alpha z^\beta \bar{z}^\gamma r^{\beta+\gamma+2\delta}\int_0^{2\pi} f(re^{i\theta}) e^{i(\gamma-\beta)\theta}   \, d\theta.
 \end{equation}
 By equation \eqref{eq:moment2}, if $r\leq \epsilon$  this last quantity equals $0$ when $\gamma-\beta\geq 1$. It follows that $Q_n(z,\bar{z})$ equals a sum that involves only terms of the form $(z\bar{z})^\alpha z^\beta \bar{z}^\gamma $ with $\gamma\leq \beta$. Hence, $Q_n$ is a holomorphic function of $z,z\bar{z}$. Because $\{Q_n\}$ converges uniformly to $f$ on $\{z\colon |z|\leq \epsilon/2\}$, taking the partial sums of the Taylor series of $Q_n$ about the origin gives the desired polynomial approximation of $f$.  (To get a holomorphic polynomial on $\C^2(z,w)$, replace $(z,z\bar{z})$ by $(z,w)$.)


\section{Flat hyperbolic or parabolic Bishop surfaces} \label{section:HyperbolicBishop}

The results in this section concern flat Bishop surfaces that are either parabolic models or hyperbolic. 



\begin{thm}\label{thm:hyperpara}
Fix $\lambda\in[1/2,\infty]$ and $\ell \geq 3$ (possibly $\ell=\infty$ or $\ell=\omega$).
Let $E(z,\bar{z})$ be $C^\ell$, real-valued, and $o(|z|^2)$. For $\lambda \neq 1/2$ define $\rho(z,\bar{z})=z\bar{z}+\lambda(z^2+\bar{z}^2)+E(z,\bar{z}),$  where $\lambda=\infty$ is interpreted as $\rho(z,\bar{z})=z^2+\bar{z}^2+E(z,\bar{z})$. If $\lambda=1/2$ define $\rho(z,\bar{z})=z\bar{z}+\frac{1}{2}(z^2+\bar{z}^2).$
For $\delta_1>0$ and $\delta_2>0$ sufficiently small, define $M=\{(z,w)\in \C^2\colon w=\rho(z,\bar{z}), |z|<\delta_1, |w|<\delta_2\}$.

    \begin{enumerate}[(i)]

\item $\CR^0(M)=C^0(M)$ has the fixed-neighborhood
approximation property, so $\CR^k(M)=\CR_P^k(M)$ for every $k\leq \ell$.

\item $\CR^\ell_P(M)\supsetneq\CR_H(M)$.

\item $\CR_H(M)$ does not have the fixed-neighborhood extension property at the origin.

\end{enumerate}
\end{thm}

\begin{proof} 
If $\lambda > 1/2$ (the hyperbolic case), the origin is the only CR singular point of $M$. If $\lambda = 1/2$ (the parabolic case), the CR singular points have the form $(it,0)$ for $t$ real. Also, $M$ is totally real away from the CR singular points.
It follows that $C^k(M)=\CR^k(M)$ for every $k\leq \ell$.
 
Note that $\CR_H(M)$ does not have the fixed-neighborhood extension property at the origin because $M$ is contained in the Levi-flat hypersurface
given by $\Im w = 0$.

We prove that $C^0(M)$ has the fixed-neighborhood
approximation property at the origin. (From this it follows that $\CR^k(M)=\CR_P^k(M)$ for every $k\leq \ell$.) Throughout we write $s=\Re w$. First we claim that,  for $\epsilon>0$ small,  the level sets of $\rho(z,\bar{z})$ in $|z|\leq \epsilon$ have a connected complement in $\C$ and empty interior. This is clear if $\lambda=1/2$, and if $\lambda\neq 1/2$ the claim follows from the Morse lemma because the Hessian of $\rho$ is nondegenerate at the origin.
Now put $K=\{(z,s)\colon |z|\leq \delta_1/2, |s| \leq \delta_2/2\}$, a compact neighborhood of the origin in $\C\times\R$.
Fix $s\in\R$ for which $K_s=\{z\in\C\colon (z,s)\in M\cap K\}$ is nonempty. Then by the claim $\C\setminus K_s$ is connected, and $K_s$ has empty interior. Thus, if $f\in C^0(M)$ we can use Mergelyan's theorem to uniformly approximate $K_s\ni z \mapsto f(z,s)$ by polynomials in $z$. This holds whenever $K_s$ is nonempty, so we may use Theorem \ref{thm:graph}
(with Remark~\ref{rmk:graph}) to uniformly approximate $f$ on $M\cap K$ by polynomials in $(z,s)$. To get holomorphic polynomials on $\C^2(z,w)$, replace $(z,s)$ by $(z,w)$. This proves that $C^0(M)$ has the fixed-neighborhood
approximation property at the origin. 

Now we show that $\CR_P^\ell(M)\neq \CR_H(M)$. 
The function $f(z,w)=\bar{z}$, considered as a function on $M$, belongs to $C^\ell(M)=\CR_P^\ell(M)$. Assume for a contradiction that $f$ can be extended to a neighborhood of the origin in $\C^2$ as a holomorphic function $g$. Then, if $L$ is the vector field on $\C^2$ defined by $L=\partial/\partial\bar{z}$, we have $Lg\equiv0$ near the origin. But $L$ is tangent to $M$ at the origin, and $L_0f\neq 0$. This is a contradiction.
\end{proof}

\begin{remark}
By the fixed-neighborhood approximation property for $C^0(M)$
where  $M$ is given by $w = z^2 + \bar{z}^2$,
every continuous function on a compact subset of $\C$ can be uniformly
approximated by polynomials in $z$ and $\bar{z}^2$.
This is a special case of a result due to Minsker~\cite{Minsker}
and later generalized by Mondal~\cite{Mondal},
whose work we mentioned in the introduction.
\end{remark}


\section{A manifold with a large hull} \label{section:Hull}

In this section we study properties of the submanifold of $\C^3$ defined
by $w=\bar{z}_1 z_2$. Note that this submanifold is locally a diffeomorphic image of
$\R^2 \times \C$ under a CR map (e.g., \cite{LNR:Severi}*{Remark 1.3}) and
is Levi-flat at CR points, so it is perhaps surprising that it has a large
disc hull and in fact satisfies the strong iterated SADH
condition at the origin (a CR singular point).

\begin{thm}\label{thm:tar}
Define  $M \subset \C^3$ by
$M=\{(z_1,z_2,w)\colon w=\bar{z}_1 z_2 ,
\snorm{z}^4+\sabs{w}^2 < \delta \}$
for a given $\delta > 0$.

\begin{enumerate}[(i)]
\item
$M$ satisfies the strong iterated SADH condition
(and hence also the DH condition)
at the origin,
and therefore at all CR singularities.
\item
$\CR^k_P(M)=\CR_H(M)$ for all $k$.
\item
$\CR_H(M)$ (and therefore also $\CR^k_P(M)$ for all $k$) has the
fixed-neighborhood extension property at the origin, 
and hence the fixed-neighborhood approximation property at the origin.
\item
For every $k$, $\CR^k(M) \supsetneq \CR^k_P(M)$, and hence $\CR^k(M)$ has neither a polynomial approximation nor an extension property at the origin (fixed-neighborhood or otherwise).
\end{enumerate}
\end{thm}

\begin{proof}
The submanifold $M$ is bounded weighted homogeneous (use $\alpha=(1,1,2)$). By rescaling it is sufficient to prove
the result for any particular $\delta$.  Therefore, suppose that
$\delta$ is large enough so that the polydisc $\Delta_C$
which we define below fits within 
the set given by
$\snorm{z}^4+\sabs{w}^2 < \delta$.

The  set of CR singular points of $M$ is $\{(z_1,z_2,w)\in M\colon z_2=0\}$. For every CR singular point $q$ of the unbounded submanifold defined by $w=\bar{z}_1 z_2$, there exists an affine biholomorphic map of $\C^3$ onto itself that sends this submanifold onto itself and the origin to $q$. Thus, the existence of an iterated shrieking disc hull neighborhood for $M$ at a given CR singular point follows from the existence of such a neighborhood at the origin. Also, note that through every CR point of $M$ there is a connected nonsingular complex curve $\Upsilon$ (a subset of a complex line) such that $\Upsilon \subset M_{CR}$ and the closure of $\Upsilon$ contains a CR singular point.

The main difficulty in the proof of the theorem is to show that
$M$ satisfies the strong iterated SADH condition at the origin.
Given this result, here is the proof of the rest of the theorem:  Parts (ii) and (iii) follow from the observations in the preceding paragraph and Corollary \ref{thm:HangesTreves}.
 Also, the function $f(z_1,z_2,w)=\bar{z}_1$, considered as a function on $M$, belongs to $\CR^\omega(M)$ (e.g., $f=w/z_2$ on $M_{\CR}$), but $f\not\in\CR_H(M)$ (otherwise, the unique holomorphic extension would equal $w/z_2$ on an open set, an impossibility). In fact, it is easy to see directly that $f\not\in \CR^\omega_P(M)$: $f$ cannot be written as a uniform limit of holomorphic polynomials on $\{(z_1,0,0)\colon |z_1|\leq \epsilon \}\subset M$. Thus, for every $k$, $\CR^k(M) \neq \CR^k_P(M)$. This proves (iv).

We now show that $M$ satisfies the strong iterated SADH condition
at the origin. 
Fix $C > 3$, write $\Delta_C$ for the closed polydisc
$\{(\xi_1,\xi_2,\omega)\in\C^3\colon \sabs{\xi_1},\sabs{\xi_2},\sabs{\omega}\leq C\}$,
and define $A_0=\Delta_C\cap M$.
It suffices to prove that $A_0$ satisfies the
the strong iterated SADH condition at the origin, which will follow from Lemma \ref{lemma:rein}.
We thus simply need to construct the iterated disc hull
$\DH^k(M)$ and show that some iteration contains
a nonempty Reinhardt domain.

{\bf First step}: We attach discs to $A_0$.
 Define 
 \begin{equation}
A_1=\Delta_C\cap \{(\xi_1,\xi_2,\omega)\colon \Im (\omega \xi_1\bar{\xi_2})=0, \; \Re (\omega \xi_1\bar{\xi_2})\geq |\xi_1\xi_2|^2,
|\xi_2|/C \leq |\xi_1| \leq C |\xi_2|\}.
\end{equation}
We prove that if $p\in A_1$ then
there exists an analytic disc $\varphi$ attached to $A_0$
through $p$.  If also $p \in A_0$, we can use a constant disc, so we
assume $p \not\in A_0$.
First fix a point 
$p=(z_1,z_2,w)\in A_1 \setminus A_0$ with $z_2\neq 0$. Then $z_1\neq 0$, so also $w\neq 0$.
 Let  $\lambda>0$ satisfy $\lambda^2=wz_1/z_2$, so $|z_1|\leq \lambda$. Define  $\varphi(\zeta)=(\lambda\zeta, w\zeta/\lambda,w)$. Note that $z_1/\lambda\in{\overline\D}$ and $\varphi(z_1/\lambda)=p$.
Clearly $\varphi$ is attached to $M$. 
 It follows from $|w|\leq C$ and
$|z_2|/C\leq |z_1|\leq C|z_2|$ that $\varphi$ is in fact attached to $A_0$: $\lambda=|wz_1/z_2|^{1/2}\leq C$, and $|w/\lambda|=|w||wz_1/z_2|^{-1/2}\leq C$.
Next, if $p=(z_1,0,w)\in A_1\setminus A_0$, then $z_1=0$, and $p$ belongs to a disc attached to $A_0$: Define $\varphi(\zeta)=(\bar{w}\zeta, \zeta,w)$.
Thus, $A_1\subset\DH(A_0)$.  This concludes the first step.

     {\bf Second step}: We attach discs to the set $A_1$ from the first step. We prove that, for $\epsilon>0$ sufficiently small, if $p=(z_1,z_2,w)$ belongs to the set  
    \begin{equation} A_2=\{(\xi_1,\xi_2,\omega)\colon |\xi_1|\leq\epsilon, \, \frac{1}{K_1} \leq |\xi_2|\leq\frac{1}{K_2}|\omega|, \, \frac{K_3}{C}\leq |\omega|\leq C\},\end{equation}
     then $p$ belongs to a disc attached to $A_1$. 
Here $2C-K_1$ and $K_2-18$ are small positive numbers that depend on $\epsilon$, and $CK_2/K_1<K_3<C^2$.

  For $\zeta\in\C$ define \begin{equation}f(\zeta,1/\zeta)=\frac{4}{9}(\zeta-1/2)(1/\zeta-1/2).\end{equation} It is easy to see that if $|\zeta|=1$ then $f(\zeta,1/\zeta)$ is real and  $1/9 \leq f(\zeta,1/\zeta)\leq 1$. 
    For the moment fix  $\lambda\in (0,1/9]$ and $\theta\in\R$. 
  Consider the map
\begin{equation} \varphi(\zeta)=(e^{i\theta}\zeta f(\zeta,1/{\zeta}),  \; w\lambda e^{i\theta} \zeta, w).\end{equation}
If $\lambda \geq 1/(C|w|)$, then $\varphi$ is attached to $A_1$: When $|\zeta|=1$ we have 
$f(\zeta,\bar{\zeta}) \geq \lambda f^2(\zeta,\bar{\zeta})>0$ because
  $0<f(\zeta,\bar{\zeta}) \leq 1\leq 1/\lambda$. Also, if $|\zeta|=1$ then
\begin{equation}
 |w\lambda e^{i\theta} \zeta|/C\leq |e^{i\theta}\zeta f(\zeta,1/{\zeta})|\leq C |w\lambda e^{i\theta} \zeta|     
\end{equation} because $\lambda \leq 1/9$, $|w|\leq C$, $1/9 \leq f(\zeta,1/\zeta)\leq 1$, and $C|w|\lambda\geq 1$. 
  
  Now we show that $\varphi(\zeta)=p$ for some $\zeta$ with $|\zeta|\leq 1$ and for some choice of $\lambda$ and $\theta$. To satisfy the first component of this equation, we use the fact that $\zeta f(\zeta,1/{\zeta})$ maps a neighborhood of $\zeta=1/2$ onto a neighborhood of $0$. In fact, from the first component we find the requirement $e^{i\theta}\zeta f(\zeta,1/{\zeta})=z_1$, 
and solving the resulting quadratic equation in $\zeta$ gives 
\begin{equation}
    \zeta=\frac{5}{4}-\sqrt{\frac{9}{16}-\frac{9}{2}e^{-i\theta}z_1}.
\end{equation}
Here the square root is chosen so that $z_1=0$ corresponds to $\zeta=1/2$. The second component of the equation  $\varphi(\zeta)=p$ then requires that 
\begin{equation} \label{eq:comptwo}
   z_2= w \lambda e^{i\theta} \left(\frac{5}{4}-\sqrt{\frac{9}{16}-\frac{9}{2}e^{-i\theta}z_1}\right).
\end{equation}
The modulus of $z_1$ is small, so on  the right side of equation \eqref{eq:comptwo}
the argument of the factor in parentheses is near 0.
We choose $\theta$ so that the right side of equation \eqref{eq:comptwo} has the same argument as $z_2$. Then we choose $\lambda\geq 1/(C|w|)$ 
so that the right side has the same modulus as $z_2$. 
Thus, $A_2\subset\DH(A_1)$. This concludes the second step.

From the first two steps we conclude that $A_2\subset\DH^2(A_0)$. Note that
is $A_0$ is bounded weighted homogeneous and that
$A_2$ contains a nonempty Reinhardt domain. By Lemma \ref{lemma:rein}, $A_0$
satisfies the strong iterated SADH condition at $0$, as desired.
%
%
%
%
\end{proof}


\section{Fixed-neighborhood
approximation without fixed-neighborhood extension} \label{section:BTnofixednbhd}

In this section we consider the submanifold of $\C^3$ defined by
$w=\sabs{z_1}^2-\sabs{z_2}^2$. The proof that $\CR^\omega$ has the
fixed-neighborhood approximation property at the origin depends on
constructing analytic discs, but the submanifold satisfies
neither
the iterated SADH condition nor the DH condition:
Because the submanifold is contained in $\C^2\times \R$, every attached disc must be contained in $\C^2\times \R$, and the same will
be true for the $\SADH_0$ hull.

\begin{thm}\label{thm:anote}
\pagebreak[2]
Define  $M \subset \C^3$ by $M=\{(z_1,z_2,w)\colon w=\sabs{z_1}^2-\sabs{z_2}^2, \snorm{z}^4+\sabs{w}^2 < \delta\}$
for a given $\delta > 0$.

\begin{enumerate}[(i)] 
\item
In the real-analytic category:
\begin{enumerate}[(a)]
\item $\CR^\omega(M)=\CR^\omega_P(M)=\CR_H(M)$
\item
$\CR^\omega(M)$ has the extension property
and the fixed-neighborhood approximation property at the origin.  
\item
$\CR_H(M)$ does not have the fixed-neighborhood 
extension property at the origin.  
\end{enumerate}

\item
In the smooth category:
\begin{enumerate}[(a)]
\item
$\CR^\infty(M)\supsetneq\CR^\infty_P(M)\supsetneq\CR_H(M)$
\item
$\CR^\infty_P(M)$ does not have the fixed-neighborhood approximation property at the origin.  
\end{enumerate}

\end{enumerate}
\end{thm}

\begin{proof}
Note that the origin is the only CR singular point of $M$.
Furthermore, the manifold $M$ is bounded weighted homogeneous, and
we wish to work in some
neighborhood of the origin.
After rescaling,
we will assume that $\delta$ is large enough.
The construction
below gives discs that can all fit within a neighborhood for a large enough
$\delta$, and we
will, for simplicity, avoid mentioning $\delta$ explicitly.

First we consider functions that are real-analytic on $M$. Fix $f\in \CR^\omega(M)$.
 Because the quadratic $|z_1|^2-|z_2|^2$ satisfies the rank condition in \cite{LNR:Severi}*{Theorem 1.1}, $f$ extends as a holomorphic function to a neighborhood of the origin in $\C^3$. Because the origin is the only CR singular point of $M$, $f\in \CR_H(M)$. Thus, 
$\CR^\omega(M)=\CR^\omega_P(M)=\CR_H(M)$, and 
$\CR^\omega(M)$ has the extension property at the origin.
Because $M$ is contained in the Levi-flat hypersurface
given by $\Im w = 0$, $\CR_H(M)$ does not have the fixed-neighborhood 
extension property at the origin.  
It remains to show that $\CR^\omega(M)$ has the fixed-neighborhood approximation property at the origin. 

We construct analytic discs. Let $A_1=\{(z_1,z_2,s)\in \C^2\times\R\colon s\geq 0,|z_1|^2\leq|z_2|^2+s\}$.
Assume that $p=(z_1,z_2,s)\in A_1$. If $z_2=s=0$, then $z_1=0$, so $p\in M$. Now assume that $|z_2|^2+s> 0$.
Then $p$ belongs to a disc attached to $M$:  Define $\varphi(\zeta)=(\zeta,z_2,s)$ for $|\zeta|^2\leq |z_2|^2+s$. Then $\varphi$ is attached to $M$. Also, $p$ belongs to the disc because $|z_1|^2\leq|z_2|^2+s$.

Similarly, if $A_2=\{(z_1,z_2,s)\in \C^2\times\R\colon s\geq 0,|z_2|^2\leq|z_1|^2-s\}$
and $p=(z_1,z_2,s)\in A_2$  with  $|z_1|^2-s>0$,
then $p$ belongs to a disc attached to $M$:  Define $\varphi(\zeta)=(z_1,\zeta,s)$ for $|\zeta|^2\leq |z_1|^2-s$. Then $\varphi$ is attached to $M$. Also, $p$ belongs to the disc because $|z_2|^2\leq|z_1|^2-s$.

 Thus, for each fixed $s_0\geq 0$, discs attached to $M$ cover $\{(z_1,z_2,s)\colon s=s_0\}$.
 The same result is true if $s_0< 0$.
 Given a disc
 $\varphi(\zeta) = (z_1(\zeta),
 z_2(\zeta),w(\zeta))$, the family
 of discs
 $\varphi_t(\zeta) = (t z_1(\zeta),
 tz_2(\zeta),t^2w(\zeta))$ will
 stay attached to $M$ and shrink to the origin
 as $t \to 0$.
 This concludes the construction of analytic discs.

Now let $f\in \CR^\omega(M)$. As we proved earlier, $f\in \CR_H(M)$.
By Lemma~\ref{lemma:simplyconn}, we can  extend $f$ as a holomorphic
function to a neighborhood of $M$ in $\C^3$. Now apply
Corollary~\ref{cor:prophullSADHflat} to extend $f$ to a fixed neighborhood (independent of $f$) of the origin in $\C^2 \times \R$,
that is, extend $f$ to a real-analytic CR function on this neighborhood. Therefore, $f$ is holomorphic in $z$ for any fixed $\Re w$, and we can make this neighborhood have the form $V \times I$ where $I$ is an interval and $V$ is a polydisc.
Write $s=\Re w$.
For fixed $s$, use the partial sums of the Taylor series of the extension to get an  approximation of $f$ by a holomorphic polynomial whose coefficients depend on $s$. Now use Theorem \ref{thm:graph}, with Remark~\ref{rmk:graph},
to get polynomials (in $(z_1,z_2,s)$, hence in $(z_1,z_2,w)$) approximating $f$ on a fixed neighborhood of the origin in~$M$. 

Now we consider the class $C^\infty(M)$. We continue to write $s=\Re w$.
First we show that $\CR^\infty(M) \neq \CR^\infty_P(M)$.
Define $f \colon M \to \C$ by
\[
f(z_1,z_2,s) =
\begin{cases}
\frac{1}{z_1} e^{-1/s^2} & \text{if $s > 0$,} \\
0 & \text{if $s = 0$,} \\
\frac{1}{z_2} e^{-1/s^2} & \text{if $s < 0$.}
\end{cases}
\] In 
\cite{crext3}*{Example 2.4}, 
it is shown that $f\in\CR^\infty(M)$ and that there is no neighborhood of the
origin in $\C^2 \times \R$ to which $f$ extends as a CR function of any regularity. It follows that $f\not\in\CR_P^\infty(M)$: If $f$ were a uniform limit on a compact neighborhood of the origin of a sequence of holomorphic polynomials, using the above construction of analytic discs attached to $M$
and adapting the proof of Theorem~\ref{thm:prophullk} would give a continuous CR extension of $f$ to a neighborhood of the origin in $\C^2 \times \R$.

Next we show that  $\CR_P^\infty(M)$ does not have the fixed-neighborhood
approximation property at the origin.
For $\epsilon>0$ define $f_\epsilon \colon M \to \C$ by
\[
f_\epsilon(z_1,z_2,s) =
\begin{cases}
\frac{1}{z_1} e^{-1/(s-\epsilon)^2} & \text{if $s > \epsilon$,} \\
0 & \text{if $s \in[-\epsilon,\epsilon]$,} \\
\frac{1}{z_2} e^{-1/(s+\epsilon)^2} & \text{if $s < -\epsilon$.}
\end{cases}
\]
Then  $f_\epsilon\in\CR^\infty(M)$. In fact, $f_\epsilon\in \CR_P^\infty(M)$:  $f_\epsilon$  is identically $0$ near the origin (the only CR singularity of $M$), and at CR points 
we apply the Baouendi--Tr\`eves approximation theorem.

Now assume for a contradiction that $\CR_P^\infty(M)$ has the fixed-neighborhood approximation property at the origin, and let $K$ be an associated compact neighborhood of the origin in $M$. 
We use the functions $f$, $f_\epsilon$ defined above.
Because $f_{1/n}\rightarrow f$ uniformly on $K$ and each $f_{1/n}$ is supposed to be a uniform limit on $K$ of a sequence of holomorphic polynomials, it  follows that  $f$ is a uniform limit on $K$
of a sequence of holomorphic polynomials. This contradiction proves that 
$\CR_P^\infty(M)$ does not have the fixed-neighborhood
approximation property at the origin.

Note also that $f_\epsilon\in\CR_P^\infty(M)\setminus \CR_H(M)$.
\end{proof}


%


\def\MR#1{\relax\ifhmode\unskip\spacefactor3000 \space\fi%
  \href{http://mathscinet.ams.org/mathscinet-getitem?mr=#1}{MR#1}}



\begin{bibdiv}
\begin{biblist}


\bib{BER:book}{book}{
  author={Baouendi, M. Salah},
  author={Ebenfelt, Peter},
  author={Rothschild, Linda Preiss},
  title={Real submanifolds in complex space and their mappings},
  series={Princeton Mathematical Series},
  volume={47},
  publisher={Princeton University Press, Princeton, NJ},
  date={1999},
  pages={xii+404},
  isbn={0-691-00498-6},
  review={\MR{1668103}},
}

\bib{B-T}{article}{
   author={Baouendi, M. S.},
   author={Tr\`eves, F.},
   title={A property of the functions and distributions annihilated by a
   locally integrable system of complex vector fields},
   journal={Ann. of Math. (2)},
   volume={113},
   date={1981},
   number={2},
   pages={387--421},
   issn={0003-486X},
   review={\MR{0607899}},
}











%

\bib{Dwilewicz}{article}{
   author={Dwilewicz, Roman J.},
   title={Global holomorphic approximations of Cauchy-Riemann functions},
   conference={
      title={Complex analysis and dynamical systems IV. Part 1},
   },
   book={
      series={Contemp. Math.},
      volume={553},
      publisher={Amer. Math. Soc., Providence, RI},
   },
   isbn={978-0-8218-5196-8},
   date={2011},
   pages={31--44},
   review={\MR{2868586}},
}


\bib{HangesTreves}{article}{
   author={Hanges, Nicholas},
   author={Tr\`eves, Fran\c{c}ois},
   title={Propagation of holomorphic extendability of CR functions},
   journal={Math. Ann.},
   volume={263},
   date={1983},
   number={2},
   pages={157--177},
   issn={0025-5831},
   review={\MR{0698000}},
}

\bib{HarveyWells71}{article}{
   author={Harvey, F. Reese},
   author={Wells, R. O., Jr.},
   title={Holomorphic approximation on totally real submanifolds of a
   complex manifold},
   journal={Bull. Amer. Math. Soc.},
   volume={77},
   date={1971},
   pages={824--828},
   issn={0002-9904},
   review={\MR{0289809}},
}

\bib{HarveyWells72}{article}{
   author={Harvey, F. Reese},
   author={Wells, R. O., Jr.},
   title={Holomorphic approximation and hyperfunction theory on a $C\sp{1}$
   totally real submanifold of a complex manifold},
   journal={Math. Ann.},
   volume={197},
   date={1972},
   pages={287--318},
   issn={0025-5831},
   review={\MR{0310278}},
}

\bib{HormanderWermer}{article}{
   author={H\"{o}rmander, L.},
   author={Wermer, J.},
   title={Uniform approximation on compact sets in $C\sp{n}$},
   journal={Math. Scand.},
   volume={23},
   date={1968},
   pages={5--21 (1969)},
   issn={0025-5521},
   review={\MR{0254275}},
}

\bib{Ivashkovich}{article}{
   author={Ivashkovich, Sergei},
   title={Discrete and continuous versions of the continuity principle},
   journal={J. Geom. Anal.},
   volume={32},
   date={2022},
   number={8},
   pages={Paper No. 226, 27},
   issn={1050-6926},
   review={\MR{4443556}},
}

\bib{Kohn-Rossi}{article}{
   author={Kohn, J. J.},
   author={Rossi, Hugo},
   title={On the extension of holomorphic functions from the boundary of a
   complex manifold},
   journal={Ann. of Math. (2)},
   volume={81},
   date={1965},
   pages={451--472},
   issn={0003-486X},
   review={\MR{0177135}},
}


\bib{crext1}{article}{
   author={Lebl, Ji\v r\'\i },
   author={Noell, Alan},
   author={Ravisankar, Sivaguru},
   title={Extension of CR functions from boundaries in ${\mathbb C}^n\times{\mathbb R}$},
   journal={Indiana Univ.\ Math.\ J.},
   volume={66},
   date={2017},
   number={3},
   pages={901--925},
   issn={0022-2518},
   review={\MR{3663330}},
}

\bib{crext2}{article}{
   author={Lebl, Ji\v r\'\i },
   author={Noell, Alan},
   author={Ravisankar, Sivaguru},
   title={Codimension two CR singular submanifolds and extensions of CR
   functions},
   journal={J.\ Geom.\ Anal.},
   volume={27},
   date={2017},
   number={3},
   pages={2453--2471},
   issn={1050-6926},
   review={\MR{3667437}},
}

\bib{crext3}{article}{
   author={Lebl, Ji\v r\'\i },
   author={Noell, Alan},
   author={Ravisankar, Sivaguru},
   title={On Lewy extension for smooth hypersurfaces in $\C^n\times\R$},
   journal={Trans.\ Amer.\ Math.\ Soc.},
   volume={371},
   year={2019},
   pages={6581--6603},
   review={\MR{3937338}}
}




\bib{LNR:Severi}{article}{
   author={Lebl, Ji\v r\'\i },
   author={Noell, Alan},
   author={Ravisankar, Sivaguru},
   title={A CR singular analogue of Severi's theorem},
   journal={Math. Z.},
   volume={299},
   date={2021},
   number={3-4},
   pages={1607--1629},
   issn={0025-5874},
   review={\MR{4329261}},
}



\bib{LNR:Cartan}{article}{
    author={Lebl, Ji\v r\'\i },
   author={Noell, Alan},
   author={Ravisankar, Sivaguru},
   title={Cartan uniqueness theorem on nonopen sets},
   journal={Complex Anal. Synerg.},
   volume={8},
   date={2022},
   number={3},
   pages={Paper No. 17, 7 pp.},
   issn={2524-7581},
   review={\MR{4476930}},
 note={\href{https://arxiv.org/abs/2112.07585v3}{arXiv:2112.07585v3} for erratum}
}


\bib{Lewy}{article}{
   author={Lewy, Hans},
   title={On the local character of the solutions of an atypical linear
   differential equation in three variables and a related theorem for
   regular functions of two complex variables},
   journal={Ann. of Math. (2)},
   volume={64},
   date={1956},
   pages={514--522},
   issn={0003-486X},
   review={\MR{0081952}},
}

\bib{Mergelyan}{article}{
   author={Mergelyan, S. N.},
   title={Uniform approximations to functions of a complex variable},
   journal={Amer. Math. Soc. Translation},
   volume={1954},
   date={1954},
   number={101},
   pages={99},
   review={\MR{0060015}},
}

\bib{Minsker}{article}{
   author={Minsker, Steven},
   title={Some applications of the Stone-Weierstrass theorem to planar
   rational approximation},
   journal={Proc. Amer. Math. Soc.},
   volume={58},
   date={1976},
   pages={94--96},
   issn={0002-9939},
   review={\MR{0467322}},
}


\bib{Mondal}{article}{
   author={Mondal, Golam Mostafa},
   title={Polynomial convexity and polynomial approximations of certain sets
   in $\mathbb{C}^{2n}$ with non-isolated CR-singularities},
   journal={J. Geom. Anal.},
   volume={33},
   date={2023},
   number={8},
   pages={Paper No. 251, 34},
   issn={1050-6926},
   review={\MR{4592426}},
}

    

\bib{NacPor}{article}{
   author={Nacinovich, Mauro},
   author={Porten, Egmont},
   title={Locally approximable CR functions, a sharp maximum modulus principle and holomorphic extension},
   note={Preprint \href{https://arxiv.org/abs/2402.19057}{arXiv:2402.19057}},
}

\bib{Severi:31}{article}{
   author={Severi, F.},
   title={Risoluzione generale del problema di Dirichlet per le funzioni biarmoniche},
   journal={Atti Accad. Naz. Lincei, Rend., VI. Ser.},
   volume={13},
   year={1931},
   pages={795--804}
}

\bib{Tumanov}{article}{
   author={Tumanov, A. E.},
   title={Extension of CR-functions into a wedge from a manifold of finite
   type},
   language={Russian},
   journal={Mat. Sb. (N.S.)},
   volume={136(178)},
   date={1988},
   number={1},
   pages={128--139},
   issn={0368-8666},
   translation={
      journal={Math. USSR-Sb.},
      volume={64},
      date={1989},
      number={1},
      pages={129--140},
      issn={0025-5734},
   },
   review={\MR{0945904}},
}










\end{biblist}
\end{bibdiv}
\end{document}